\documentclass[a4paper,10pt]{article}

\def \Z {{\mathbf {Z}}}
\def \R {{\mathbf {R}}}
\def \N {{\mathbf {N}}}
\def \F {{\bf F}}
\def \B {{\cal B}}
\def \T {{\hat T}}
\def \Q {{\bf Q}}

\def \V {{\bf V}}
\def \P {{\bf P}}

\def \tt {{\hat {\bf T}}}
\def\uu{\bigsqcup}
\def\eps{\varepsilon}
\title{ О сохраняющих меру преобразованиях ранга  один }
\author{В.В. Рыжиков}
\date{25.05.2020}
\usepackage[T1]{fontenc} 
\usepackage[cp1251]{inputenc} 
\usepackage[russian]{babel}

\begin{document}

\maketitle
\begin{abstract}{   Преобразования ранга один служат источником примеров в эргодической теории, показывающих
разнообразие алгебраических, асимптотических и  спектральных свойств динамических систем. 
Свойства преобразований ранга один   тесно связаны со структурой полугруппы слабых пределов   его степеней.   
В этом ключе  изучаются     известные и новые конструкции преобразований. }
\end{abstract}
\large

\section{Введение} 

Преобразованием в статье  называется  сохраняющее меру обратимое 
преобразование   пространства Лебега $(X,\B, \mu)$. Пока речь идет о 
стандартном вероятностном пространстве, позже будут также рассмотрены преобразования 
пространства  с бесконечной мерой.
Преобразования образуют группу Aut, на которой естественно возникает полная метрика. 
В.А. Рохлин и П. Халмош доказали типичность в смысле Бэра слабого перемешивания
при отсутствии сильного. Тот факт, что периодически преобразовния плотны   в Aut, 
привел к методу  аппроксимаций периодическими преобразованиями,  развитому 
А.Б.Катком, В.И.Оселедцем, А.М.Степиным и др. 
для обнаружения типичных свойств преобразований и построения  примеров.  
Преобразование ранга один допускают аппроксимацию в смысле работы \cite{KS}, но  его определение   
 свободно от привлечения внешних,  аппроксимирующих преобразований.
 
Преобразование  $T$ обладает \it  рангом один, \rm если  некоторая последовательность измеримых 
разбиений фазового вида
$$\xi_j=\{E_j, TE_j, T^2E_j,\dots, T^{h_j-1}E_j, \tilde{E_j}\},$$
стремится к разбиению на точки.  Это означает, что любое  множество 
$A\in\B$ приближается некоторой последовательностью  $\xi_j$-измеримых множеств  $A_j$: 
  $$\mu(A\Delta A_j)\to 0, \ j\to\infty.$$\rm
Последовательность разбиений можно 
модифицировать так, что полученная новая последовательность будет   монотонной \cite{Bax}. 
Это приводит к другому определению ранга один, в котором преобразование является конструкцией,
однозначно определенной заданными параметрами.

\vspace{3mm}
\bf Конструкция преобразования ранга один. \rm   
Пусть дана последовательности целочисленных векторов 
$$ \bar s_j=(s_j(1), s_j(2),\dots, s_j(r_j-1),s_j(r_j)), \ r_j>1, $$

 Положим $h_1=1$.  Мы будем индуктивно доопределять фазовое пространство и преобразование.
На каждом новом  этапе то, что было определено на предыдущих, не меняются в дальнейшем.
 
На этапе $j$  частично определенное  преобразование    $T$ является обычной
 перестановкой непересекающихся интервалов, образуя из них  
башню  вида
$$E_j, TE_j T^2,E_j,\dots, T^{h_j-1}E_j.$$
  Разрежем интервал  $E_j$ на  $r_j$ интервалов  $E_j^1, E_j^2,  E_j^3,\dots, E_j^{r_j}$ одинаковой меры.
Рассмотрим  колонны 
$$E_j^i, TE_j^i ,T^2 E_j^i,\dots, T^{h_j-1}E_j^i, \  i=1,2,\dots, r_j.$$

Над  колонной с номером $i$ надстраиваем   $s_j(i)$ новых интервалов и получим набор интервалов 
$$E_j^i, TE_j^i, T^2 E_j^i,\dots,  T^{h_j+s_j(i)-1}E_j^i$$
(интервалы  не пересекаются и имеют одинаковую меру).
Для всех   $i<r_j$ положим 
$$T^{h_j+s_j(i)}E_j^i = E_j^{i+1}.$$ 
Тем самым мы сложили колонны в  башню  этапа $j+1$
$$E_{j+1}, TE_{j+1} T^2 E_{j+1},\dots, T^{h_{j+1}}E_{j+1},$$
где   $$E_{j+1}= E^1_j,\
h_{j+1}=h_jr_j +\sum_{i=1}^{r_j}s_j(i).$$
Продолжая построение, получим на объединении $X$ всех интервалов  сохраняюще меру преобразование $T$.
Если мера  $Х$ конечна,  нормируем ее. 

Аналогичным образом определяется более широкий класс конструкций, кода интервалы 
 под действием $T$  могут менять длину.
 Такое преобразование  с квазиинвариантной мерой, не эвивалентной никакой инвариантной мере,
 предъявлено в работе \cite{O1}. 
Следующие примеры хорошо известны в эргодической теории:

(i) преобразования Чакона \cite{Ch},  \cite{Ch69} 
$$\bar s_j=(2,3,1), \ \bar s_j=(0,1);$$

(ii) стохастические конструкции Орнстейна \cite{O} 
$$s_j(i)= b_j+a_j(i)-a_j(i+1),   1\leq a_j(i)\leq b_j\to \infty;$$

(iii) преобразование  дель Джунко-Рудольфа \cite{JR}  

$$\bar s_j=(0,0,\dots, 0,0,1,0,0,\dots, 0,0), \ r_j\to\infty;$$

(iv) конструкция Катка \cite{Katok}   

$$\bar s_j=(0,0,\dots,0, 0,1,1,\dots, 1,1), \ r_j\to\infty.$$

В статье, в частности,   рассматриваются следующие свойства преобразований: 
перемешивание, слабое кратное перемешивание, тривиальность централизатора, 
отсутствие факторов, жесткость,   минимальные самоприсоединения, простота самоприсоединений, 
дизъюнктность сверточных  степеней спектральной меры, простота спектра симметрических степеней.  

Свойства конструкций полностью определяются последовательностью надстроек  $\bar s_j$,  которые 
задают крышу над башней. 
По этой причине конструирование  преобразований ранга один имеет аналогию 
с методом, описанным  в художественной литературе: \it
there was a most ingenious 
architect who had contrived    a new method for building houses, by beginning at the roof, 
  and working downward to the foundation  \rm  (Jonathan Swift, Gulliver's Travels).


\section{Обзор результатов}
\bf Перемешивающие конструкции.  \rm 
Преобразование  $T$ вероятностного пространства 
обладает свойством перемешивания, если для любых $A,B\in\B$  выполнено
$$\mu(T^iA\cap B)\to \mu(A)\mu(B),\ i\to\infty.$$  
Класс  преобразований ранга один, задаваемых  случайным параметром, 
 введен   в работе \cite{O}. Д.Орнстейн  доказал, что параметрическое  большинство      
в этом классе образуют  перемешивающие преобразования с  тривиальным централизатором. 
 Перемешивающие преобразования ранга один обладают кратным  перемешиванием  \cite{Ka},  \cite{R92}.
 
В \cite{B} предложен новый подход для изучения 
спектров преобразований ранга один  и установлена  сингулярность спектра 
 для параметрического  большинства
 орнстейновских преобразований. Общий результат о сингулярности спектра конструкций преобразований  с медленно растущей
последовательностью $r_j$ дан в \cite{KR}. 
  Свойство перемешивания для  класса лестничных конструкций  доказано в \cite{Ad}. 
Небольшие модификации лестничных конструкций нашли приложения в исследовании кратностей спектра \cite{06}.

 \bf Минимальные самоприсоединения. \rm Д.Рудольф  \cite{Ru} ввел в рассмотрение  преобразование  ранга один 
с  минимальными самоприсоединениями (MSJ). Предложенное им  свойство гарантирует 
тривиальность централизатора и отсутствие факторов (нет собственных инвариантных $\sigma$-алгебр). Следствием 
этого свойства является  минимальность  централизатора и минимальность структуры факторов  у всех  
 декартовых степеней преобразования. Минимальность,   неформально говоря,   означает присутствие  только 
тех объектов, которые заведомо существуют.  
Преобразование  $T$ с минимальными самоприсоединениями  применялось для построения разнообразных контрпримеров,
коллекция которых приведена в  \cite{Ru}.
Рассмотрим произведения $$T'=T\times T\times T\times \dots, \ \ T''= T\odot T\times T\times \dots,$$
 где $T\odot T$ -- фактор призведения $T\times T$, являющийся ограничением 
$T\times T$ на алгебру множеств, неподвижных относительно симметрии $(x,y)\to (y,x)$.  Тогда $T'$ и $T''$ слабо
изоморфны в смысле Синая:  $T'$ изоморфен фактору  $T''$ и наоборот. Однако, они  не изоморфны,
что  есть следствие свойства MSJ.   Другой пример: произведение  $T\times T\times T$  имеет корень  степени 3,
 но не имеет  корней других степеней. 

{\it Cамоприсоединением порядка 2} преобразования $T$,  называется мера $\nu$ на 
$X\times X$, для которой выполнено   $ (T\times T)\nu = \nu$ и  
  $\nu(X\times A) = \nu(A\times X) = \mu(A)$ для любого $A\in\B$.

Преобразованию $S$,  коммутирующему  с преобразованием $T$, отвечает  самоприсоединение  
$\Delta_S$ на $(X\times X)$, определенное
равенством $\Delta_S=(Id\times S)\Delta$, где $\Delta$ -- диагональная мера: 
$\Delta(A\times B)=\mu(A\cap B)$. 

Если  у преобразования нет  эргодических  самоприсоединений, отличных от $\Delta_{T^n}$ и $\mu\times\mu$,  
говорят, что оно обладает свойством MSJ(2), минимальных самоприсоединений порядка два. 
Близкое понятие простоты присоединений,  обозначаемое Simpl(2), означает, что
 нет  эргодических  самоприсоединений, отличных от $\Delta_{S}$ и $\mu\times\mu$. 

Аналогично  определяются минимальные самоприсоединения больших  порядков. 
Неформальное определение MSJ таково: преобразование $T$  допускает только очевидные    самоприсоединения 
всех порядков. Не отвлекаясь на интересные общие вопросы теории присоединений, отметим, что для групповых действий 
 MSJ(4), для преобразований  MSJ(3), а для потоков   MSJ(2) эквивалентны   MSJ. Для знакомства 
с этой проблематикой  см., например, \cite{K4},\cite{GHR},\cite{96},\cite{97}.

\bf  Типичные свойства.  \rm Типичные в Aut свойства преобразований  суть
свойства некоторых преобразований ранга один.  Сравнительно недавно обнаружилась категорная типичность   
свойства  MSJ.   С.В. Тихонов \cite{T1} переоткрыл метрику С. Альперна \cite{Al} на  пространстве  перемешивающих 
преобразований  { Mix}, доказав  полноту   Mix  и установив  типичность таких свойств как сингулярность 
спектра и кратное перемешивание. 
А.И. Баштанов   обнаружил  типичность ранга один в пространстве  { Mix} \cite{Ba} и 
с учетом  результатов работ  
\cite{Ka}, \cite{R92},\cite{K1}  типичность   свойства MSJ.

 В группе Aut типичные преобразования обладают  корнями \cite{K2}, они  являются групповыми расширениями  \cite{A1}
(неизвестно, является ли типичное преобразование  относительно слабо перемешивающим расширением \cite{U06}).
 Типичное преобразование  включается в поток \cite{LR},  включение  неоднозначно, централизатор типичного преобразования 
содержит бесконечномерные торы \cite{SE}.
 Все это  не имеет места для пространства Mix.  Здесь типичны   тривиальный централизатор, и отсутствие факторов.
 
\bf  Неперемешивающие конструкции.  \rm   
Конструкции с надстройками вида $\bar s_j=(2,3,1)$ и  $\bar s_j=(0,1)$ расматривались в \cite{Ch},  \cite{Ch69}
  как  примеры неперемешивающих преобразований ранга один без корней и с непрерывным спектром. 
Преобразование с $\bar s_j=(0,1)$  называют историческим преобразованием Чакона, позже мы обсудим
подробно его свойства.
Для модифицированного преобразованием Чакона  ( надстройки $\bar s_j=(0,1,0)$) 
установлено свойства  MSJ  \cite{JRS} и  дизъюнктности сверточных степеней спектральной меры  \cite{2000},
описана  структура слабого замыкания  \cite{JPRR}.

Естественным обобщением преобразований Чакона   являются \it ограниченные конструкции \rm
(конструкции ранга один, у которых все параметры ограничены). 
В \cite{2012} доказано, что все нежесткие слабо перемешивающие ограниченные конструкции обладают 
минимальными самоприсоединениями, а  в связи с работой \cite{B1}    установлена попарная дизъюнктность 
всех положительных степеней слабо перемешивающих ограниченных  конструкций.
Доказательство спектральной  дизъюнктности  степеней слабо перемешивающих ограниченных конструкций 
  получено в \cite{ALR}.   

Класс неперемешивающих преобразований образуют \it полуограниченные конструкции \rm 
(ограничены некоторые параметры).  Такое свойство как простота спектра 
степеней $\hat T^{\odot n}$  для многих полуограниченных конструкций устанавливается  проще, 
чем для ограниченных.
 В \cite{2019} построено полуограниченное  преобразование $T$  ранга один, 
у которого    тензорные произведения  
 $\hat T\otimes \hat T^{m}$ имеют простой сингулярный спектр при $1<m<2020$, а спектр 
произведений $\hat T\otimes \hat T^{n}$ при $n\geq 2020$ является счетнократным лебеговским, см. \cite{2019}.  
 Полуограниченное самоподобное   преобразование ранга один из \cite{2013c} индуцирует гауссовский автоморфизм 
$G$ с сингулярным  спектром такой, что   степень $G^3$ изоморфна тензорной степени  $G\otimes G\otimes G$. 
У такого автоморфизма $G$ множество
спектральных кратностей  степени $G^{3^p}$ есть  $\{3^p,\infty\}$. Но тогда  у гауссовского потока $G_t$
множество спектральных кратностей его элементов  непостоянно.

\bf  Слабое замыкание  действия ранга один. \rm
В  эргодической теории  хорошо известно следующее приложение слабых пределов степеней преобразования.
В работах  В.И.Оселедца и А.М.Степина 
(см. \cite{Os},\cite{St})  
 в связи с проблемой Колмогорова о групповом 
свойстве спектра преобразования  рассматривались 
слабые пределы  вида $aI+(1-a)\Theta$,  где $\Theta$ -- ортопроекция на пространство констант в $L_2$.

Слабые пределы  $aI+(1-a)\T$ применялись   в теории присоединений  \cite{96},\cite{97},  
позднее они нашли приложения  в спектральной теории динамических систем (см. \cite{Da}).   
В случае перемешивания, означающего сходимость $\T^i\to \Theta$, где $\T$ -- оператор, индуцированный 
преобразованием $T$,  нетривиальные слабые  пределы отсутствуют. Тем не менее можно использовать 
технику слабых пределов
 даже в случае перемешивания \cite{1999},\cite{07}.
В терминах нестандартного анализа это объясняется тем, что слабое замыкание степеней перемешивающего преобразования 
содержит  бесконечно близкие оператору  $\Theta$  пределы, которые вынуждают нужное спектральное свойство.

В случае вероятностного пространства доказательство свойства перемешивания преобразования 
ранга один является  нетривиальной задачей, она требует  контроля поведения величин
$\mu(T^iA\cap B)$  для всех больших $i$.   Более общей и  сложной   задачей в классе слабо перемешивающих действий
  является  описание полугруппы всех слабых пределов степеней преобразования,  см. \cite{JPRR},\cite{2012},\cite{2011d}. 
Если  слабое замыкание степеней слабо перемешивающего    преобразования $T$ ранга один полиномиально, т.е. 
состоит из линейных комбинаций степеней $\T$, то такое преобразование $T$ обладает свойством MSJ.
 
В случае бесконечной меры  фазового пространства, задача описания слабого замыкания действия ранга один может быть
 значительно  проще. Например, для сидоновских конструкций  свойство перемешивания ($\T^i\to_w 0$) 
и быстрое убывание коэффициентов 
Фурье спектральной меры   вытекает  из определения этих конструкций \cite{mmo14}.
Сравнительно легко описываются нетривиальные полугруппы  некоторых полуограниченных  неперемешивающих
конструкций, см. \cite{2019},\cite{KuR}. 

Понятие преобразования ранга один имеет такие обобщения как  конечный ранг, локальный ранг (см. \cite{Ro}),
групповые действия ранга один. Привлечение новых инвариантов  расширяет круг вопросов  
и возможности исследователей.    А. дель Джунко  использовал групповое действие ранга один 
с минимальными самоприсоединениями, 
для обнаружения простого   преобразования \cite{J} без минимальных факторов. 
Его метод был развит и нашел ряд  новых приложений \cite{D08}.
О.Н. Агеев также воспользовался похожим приемом,   рассмотрев типичное  действие специально подобранной группы,
в котором поддействие обладает желаемым спектральным свойством. Это привело к обнаружению
эргодических преобразований с заданной кратностью  однородного спектра \cite{Ag}.

 \bf Основное содержание статьи.   \rm В работе изучаются    связи слабого замыкания  действия  
со спектральными свойствами действия и структурой его самоприсоединений. Это направление наиболее 
актуально для преобразований ранга один.   С ростом ранга    взаимосвязи  становятся  
слабее или вовсе пропадают.   Приведем  дальнейший план работы.

(\S 3) Слабое замыкание действия и марковский централизатор. 
Если неразложимый марковский оператор коммутирует
с преобразованием  $T$ ранга  один, то  в слабой топологии он приближается, неформально говоря,  частями 
двух  степеней преобразования. 

(\S 4) Ограниченные конструкции. Рассмотрена   простейшая ограниченная конструкция с непрерывным 
спектром -- классическое преобразование Чакона. Доказан ряд  его свойств.

(\S 5)  Полуограниченные конструкции.  Модификации преобразования Чакона обнаруживают феномен 
срытой жесткости и наличия всех полиномиальных пределов. Обсуждается стохастическое преобразование Чакона. 

(\S 6) Стохастические конструкции.  Объясняется методика доказательства 
свойства перемешивания орнстейновских преобразований.

(\S 7) Типичность ранга один, нестандартное перемешивание, нетипичные преобразования.  
Показано, что класс сопряженности орстейновского  преобразования $T$ всюду 
плотен в пространстве  Mix.   При этом  идея случайной кограницы 
 используется дважды: она  заложена  в самой  орстейновской конструкции  $T$  
и виртуальным  образом  применяется при выборе  
подходящего сопряжения. Затрагивается вопрос об инвариантах в классе перемешивающих конструкций.   
Дано достаточное условие  на компакт, чтобы он не пересекался с классами сопряженности, содержащими 
плотное $G_\delta$-множество.

(\S 8) Алгебраические надстройки вместо случайных. Поля Галуа --  источник квазислучайных последовательностей,
позволяющий дерандомизировать стохастические конструкции. Свойства конструкций сходны с орнстейновскими, но
теперь алгебра заменяет статистику.

(\S 9) Лестничные  конструкции. Объясняется  методика  доказательства
свойства перемешивания  лестничных конструкций.

(\S 10)  Явная перемешивающая конструкция с двукратным  спектром.
    Известно, что существуют перемешивающие лестничные конструкции
$T$ такие, что  $T\times T$  имеет однородный спектр кратности 2. Однако, до сих пор  не был указан
конкретный  пример. Мы даем оценку  скорости роста параметров $r_j$, когда  
нужный эффект имеет место.  Например, годится  $r_j=[\ln (j +8)]$.

(\S 11) Бесконечные преобразования, самоподобные конструкции.  
 Предлагаются простые   бесконечные конструкции с необычными свойствами, например,
с самоподобием спектра.   Такие  конструкции  интересны в качестве   приложений  
к   пуассоновским и гауссовским действиям.
Дан  простой пример преобразования, не сопряженного своему обратному. Обсуждаются сидоновские конструкции,
у которых декартов квадрат диссипативен.

(\S 12) Заключительные замечания, общие и специальные 
вопросы о преобразованиях ранга один завершат статью.


\section{  Слабое замыкание действия и марковский централизатор}

\bf  Централизатор и факторы. \rm Теорема Кинга \cite{K2} утверждает, что  для преобразования $S$, коммутирующего с   
преобразованием $T$ ранга один,   для некоторой последовательности $k(j)$ верно, что   
$\T^{k(j)}\to \hat S$ (здесь речь идет о слабой сходимости,
которая в данном случае совпадает с сильной операторной сходимостью). 
Следствие: если $S$ не является степенью преобразования $T$, 
то $T$ является жестким, так как $\T^{k(j+1)-k(j)}\to I$.
Любой собственный фактор преобразования 
ранга один является жестким \cite{K2}. Теорема имеет ряд применений (см. \cite{FRW}, \cite{JRR}).
В связи с ней возникает общий вопрос:
какова связь слабого замыкания действия с коммутирующим с ним марковским оператором ?

\bf Самоприсоединения  и марковский централизатор. \rm
Мы знаем, что для преобразования ранга один $T$ и  самоприсоединения вида $\nu= \Delta_{S}$ 
найдется последовательность ${k(j)}$ такая, что $\Delta^{k(j)} \to \nu$.
 В общем случае доказано только частичное приближение:  

\it если $\nu$ -- эргодическое самоприсоединение порядка 2 преобразования $T$  ранга один, то
 для некоторой  последовательности $k(j)$ имеем
$\Delta^{k(j)} \to \frac 1 2 \nu +\dots$. \rm

Следствие: 
 \it  перемешивающие преобразования ранга один обладают свойством MSJ(2).\rm

Формулировка остается верной в случае бесконечных пространств с мерой, см. \cite{2019}.
Таким образом, в случае вероятностного пространства для 
$\Z$-действий ранга один свойство минимальных самоприсоединений эквивалентно тому, что
их слабое замыкание содержится в выпулой комбинации элементов действия и 
ортопроекции $\Theta$ на пространство  констант в $L_2(X,\B,\mu)$. Бедная структура слабого замыкания
влечет за собой бедность структуры самоприсоединений со всеми вытекающими последствиями.  Наличие нетривиальных
полиномиальных пределов приводит к ряду спектральных эффектов.

Оператор $P$  в $L_2(X,\mu)$ называем марковским, если $P$ положителен, т.е. сопоставляет 
неотрицательным функциям  неотрицательные,  и операторы  $P, P^\ast$ сохраняют интеграл. 
Самоприсоединениям преобразования $T$  отвечают марковские операторы, 
коммутирующие с $\T$. Эта связь задается формулой 
$$  (Pf,g)=\int_{X\times X}  f\otimes g d\nu. $$ 

На операторном языке упомянутое выше утверждение звучит так: 
$$T^{k(j)} \to \frac 1 2 P + P',$$
где марковский  оператор $P$ отвечает  эргодическому самоприсоединению $\nu$,  а $P'$ --  
некоторый положительный оператор (в случае бесконечного пространства он может оказаться нулевым).
 Положительность здесь означает $P'f\geq 0$  при $f\geq 0$.  Оператор $P$ неразложим:
 его нельзя представить в виде полусуммы различных марковских операторов, коммутирующих с  $\hat T$.
 Эргодическому самоприсоединению преобразования отвечат оператор, являющийся крайней  точкой
 в его марковском централизаторе.

Следующее утверждение (рассматривается случай  вероятностного  пространства)  
показывает, что эргодическое самоприсоединение приближается частями мер $\Delta^{k(j)}$. 

\it  Пусть $\nu$ -- эргодическое самоприсоединение преобразования $T$  
ранга  один.
Найдется последовательность $k_j$   такая, что для любых $A,B\in\B$
$$\nu(A\times B) =\lim_j \mu(T^{k_j}A\cap B\ | Y_j)$$
где $Y_j$ является частью башни этапа $j$, а мера $Y_j$ не меньше $\frac 1 2$.  \rm

Из доказательств \cite{R92} видно, что самоприсоединение приближается двумя частями: 
 $$\nu(A\times B) =\lim_j \left(\mu(T^{k_j}A\cap B\cap Y_j) +\mu(T^{{k'}_j}A\cap B\cap {Y'}_j)\right),$$
где $Y_j,{Y'}_j$ -- непересекающиеся части башни  этапа $j$, а сумма их мер стремится к 1. 
Приведем эквивалентную  операторную  формулировку.

\bf  Теорема 3.1. \it Пусть $T$ -- преобразование ранга один, $P$ -- неразложимый марковский оператор, 
коммутирующий с  $\hat T$. Найдутся последовательности $k_j$, $k_j'$  и последовательность 
множеств $Y_j\in \B$, для которых выполняется
$$\hat Y_j\hat T^{k_j}+(I-\hat Y_j)\hat T^{k'_j}\ \to_w \ P,$$  
где $\hat Y_j$  -- оператор умножения на индикатор множества   $Y_j$.\rm 

Доказательство.  Положим 

$$ Y^k_j=\uu_{i=0}^{h_j-k} T^{i+k} E_j, \ \ a_j^k =\frac { \nu( T^k E_j\times E_j)}{\mu(E_j)}, \ \  0\leq k < h_j,$$

$$ Y^k_j=\uu_{i=0}^{h_j+k}T^i E_j, \ \ a_j^k =\frac {\nu( E_j\times T^{-k}E_j)}{\mu(E_j)},\ \ -h_j< k<0.$$

Проверяется, что для  для самоприсоединения  $\nu$  преобразования $T$ 
и $\xi_j$-измеримых множеств $A,B\in\B$  имеет место равенство
$$\nu(A\times B) = \nu_j(A\times B):=\sum_k a^k_j \mu(T^kA\cap B\cap Y^k_j).$$

Расcмотрим нормированые меры  $\Delta^k_j$, заданные формулой   
$$\Delta^k_j (A\times B)=\mu(T^kA\cap B\ | Y^k_j)$$
   при           $|k|<(1-\delta) h_j$, где  $\delta>0$ мало.

 Если самоприсоединение  $\nu$ эргодическое, то  
для весового большинства номеров $k$ (большинство подразумевается относительно весов $a^k_j$ при фиксированном $j$)
   верно, что   $$\Delta^k_j (A\times B)\approx \nu(A\times B).$$
Причина в том, что при  близости  выпуклой суммы почти инвариантных нормированных мер  к нормированной эргодической мере, 
 большинство этих мер будут близки к эргодической мере.  Иначе несложно показать, 
что мера $\nu$ не является крайней точкой в пространстве нормированных  инвариантных мер   (см.\cite{R92}, \cite{JRR}).
Нужную последовательность $\Delta^{k_j}_j\to \nu$  всегда можно выбрать так, что 
$\frac {|k_j|} {h_j} \to b\geq \frac 1 2$.
Эта возможность обеспечивается  проекционными  свойствами меры $\nu$, которые приводят к  
$\sum_k a^k_j >a>0,$   при           $2|k| >1 +\delta $.  Если выбирать $k_j$  так, что $b$ будет  максимальным, 
в случае  $b<1$
найдется антипод   $k_j'$ такой, что    $|k_j| +|k_j'|  \approx h_j$ и выполняется  $\Delta^{k_j'}_j\to \nu$.
Это приводит (случай положительных $k_j$) к  сходимости
$$\nu(A\times B) =\lim_j \left(\mu(T^{k_j}A\cap B\cap Y_j^{k_j}) + 
\mu(T^{k'_j}A\cap B\cap Y_j^{h_j-k_j})\right).\eqno(2.1) $$

Осталось заметить, что  $Y_j^{h_j-k_j}$ можно заменить на $X\setminus Y_j^{k_j}$, и
переформулировать (2.1) в терминах операторов.    Теорема доказана.

\bf Случай плоской крыши. \rm Если   выполнено  условие 
$$\mu(T^{h_j}E_j|E_j)\to 1,$$ 
 можно выбрать  $k_j'=h_j-k_j$,  тогда (2.1) превращается в
$$\nu(A\times B) =\lim_j \mu(T^{k_j}A\cap B). $$
Этот результат известен как теорема Кинга о слабом замыкании степеней преобразования с плоской крышей, 
 см. \cite{JRR}.   


\section {Ограниченные конструкции, спектр и самоприсоединения}  

Рассмотрим простейшие конструкции ранга один с ограниченными параметрами: 
преобразование Какутани-фон Неймана (одометр) и  классическое преобразование Чакона.

\bf Одометр. \rm Пусть дана последовательность $\{r_j\}$ натуральных чисел
$r_j>1.$ Опишем конструкцию с  нулевыми векторами надстроек $\bar s_j=(0,0,\dots,0)$.
Фазовое пространство $X$ в этом примере является   интервалом единичной длины
(читателю удобно представлять его полуинтервалом). На этапе $j=1$ имеем отрезок $E_1$, $h_1=1$.

Пусть на  этапе $j$  частично определенное  преобразование    $S$ является  
перестановкой непересекающихся интервалов 
$E_j, SE_j S^2,E_j,\dots, S^{h_j-1}E_j.$ На последнем интервале оно пока не определено.
Переходим  к шагу с номером $j+1$.  Представим 
интервал  $E_j$ как объединение  $r_j$ интервалов  $E_j^1, E_j^2,  E_j^3,\dots, E_j^{r_j}$
одинаковой длины.  
Для  $i=1,2,\dots, r_j$ рассмотрим  колонны $E_j^i, SE_j^i ,S^2 E_j^i,\dots, S^{h_j-1}E_j^i$
и положим $S^{h_j}E_j^i = E_j^{i+1}$ для каждого  $i<r_j$.
Получаем башню  
$$E_{j+1}, SE_{j+1} S^2 E_{j+1},\dots, S^{h_{j+1}-1}E_{j+1},$$
где   $E_{j+1}= E^1_j,\ h_{j+1}=h_jr_j.$
Преобразование $S$ не определено только на  последнем интервале.   Продолжая описанный процесс  до бесконечности,
определяем  $S$ на всем  $X$.  

Построенное преобразование  эргодично, так как инвариантное множество $A$ положительной меры устроено на интервалах
одинаково, следовательно
имеет полную меру. Этот аргумент применим ко всем преобразованиям ранга один. 
Действительно,  некоторый интервал из башни  этапа $j$ для больших $j$
в основном  состоит из элементов  множества $A$, но в силу инвариантности $A$ сказанное верно и для других 
интервалов башни. 

 Одометр $S$  индуцирует унитарный оператор $\hat S$, $\hat Sf(x)=f(Sx)$  
в пространстве $L_2(X,\B,\mu)$.
Спектр оператора   $\hat S$ -- группой, порожденная числами вида
$e^{{2\pi mi} /{r_1r_2\dots r_n}}$.   Последовательность $\hat S^{h_j}$ сходится сильно к тождественному оператору $I$
(напомним, что такие преобразования называют жесткими).
Отметим,  что все эргодические преобразования с дискретным спектром являются жесткими и обладают рангом один \cite{J76}. 

\bf Изоморфизм степеней. \rm В случае, например, когда $r_j=3$ для всех $j$, 
 соответствующий одометр  изоморфен своему квадрату. О.Н. Агеев показал, что существуют 
преобразования $R$ ранга один с непрерывным спектром с аналогичным свойством \cite{A2}.  Явные конструкции
вряд ли легко предъявить, в силу результатов \cite{2012} их параметры будут неограниченными. 
 Как следствие  теоремы Кинга \cite{K3} получаем, что  все коммутирующие с таким $R$ преобразования  имеют общее 
сопряжение со своими квадратами. Такие $R$ не могут коммутировать с  нетождественным преобразованием четного периода.

\bf Индуцированное преобразование и самоприсоединения. \rm В препринте   \cite{CK} рассматривались  
преобразования  $T=S_A$, индуцированные некоторым одометром $S$ на специально выбранном подмножестве $A\subset X$. 
Напомним, что $S_Ax = S^{n(x)}x$, где $x\in A$, а  $n(x)$ --  время первого возвращения $x$ в $A$. Авторы
 получили пример жесткого преобразования ранга один без факторов,  
но обладающего  обширной структурой самоприсоединений,  являющихся  относительно слабо перемешивающими. 
Поясним, что это означает. Пусть  марковский оператор $P$ отвечает 
самоприсоединению $\nu\neq \Delta$,  а самоприсоединение 
$\eta$ отвечает оператору $P^\ast P$.   
 Если $\eta$ эргодично,  то говорят, что  $\nu$ является относительно слабо перемешивающим расширением 
исходной системы.   Всегда ли у типичного преобразования найдется такое самоприсоединеие?  Даже в случае,
если оно задается фактором, т. е. при дополнительном условии $P^\ast=P=P^2$, это неизвестно, см. \cite{U06}. 

\bf Классическое  преобразование Чакона.  \rm
 Пусть на  нулевом этапе дан некоторый интервал $E_0$.  
На  этапе $j$     имеется башня 
$$E_j, TE_j T^2,E_j,\dots, T^{h_j-1}E_j,$$
состоящая из $h_j$ непересекающихся  интервалов одинаковой длины.
Преобразование $T$ определено на всех интервалах, исключая последний, как обычный  перенос интервалов.
На шаге $j+1$  представим 
  $E_j$ как объединение  непересекающихся   интервалов $E_j=E_j^1$ и $E_j^2$ одинаковой длины. 
Башню этапа $j$ разрезаем на две колонны, добавляя  надстройку $T^{h_j}E_j^2$
над второй колонной:

$E_j^1, TE_j^1 ,T^2 E_j^1,\dots, T^{h_j-1}E_j^1,$ \ 

$E_j^2, TE_j^2 ,T^2 E_j^2,\dots, T^{h_j-1}E_j^2, T^{h_j}E_j^2.$ \\
Положив $T^{h_j}E_j^1 =E_j^2, \ E_{j+1}= E^1_j,$
 получим  башню этапа $j+1$  из $h_{j+1}=2h_j+1$ интервалов
$$E_{j+1}, TE_{j+1} T^2 E_{j+1},\dots, T^{h_{j+1}-1}E_{j+1}.$$

Преобразование $T$ на последнем интервале пока не определено.
Продолжая построение, мы доопределяем и  фазовое пространство $X$ и  преобразования $T$.
 Отметим, что сумма мер всех добавленных интервалов  конечна, в качестве инвариантной вероятностной меры 
 рассматривается нормированная мера Лебега на $X$.
  
Приведем список 
свойств классического преобразования Чакона:

(i) слабое перемешивание (нет собственных функций кроме констант), 
отсутствие  сильного перемешивания;

(ii) спектральная мера  $\sigma_T$ преобразования взаимно сингулярна со  
сверткой $\sigma_T\ast\sigma_T$, произведение $\hat T\otimes \hat T$ имеет однородный спектр
кратности 2;

(iii)  преобразование $T$  обладает  кратным слабым перемешиванием (ниже формулировка свойства WMix(2)): 
если $\T^{m(i)},\T^{n(i)},\T^{m(i)-n(i)}\to_w\Theta$,
то для любых $A,B,C\in\B$
$$\mu(A\cap \T^{m(i)}B\cap\T^{n(i)}C)\ \to \mu(A)\mu(B)\mu(C);$$

(iv)   преобразование $T$  обладает свойством MSJ, следовательно,  у него тривиальный централизатор и  нет факторов.

Покажем, как устанавливаются перечисленные  свойства преобразования. 

(i). Отсутствие перемешивания вытекает из наличия 
слабого предела:
 последовательность $\hat T^{-h_j}$  слабо сходится к  оператору
$$P(\T)=\sum_{k=1}^{\infty} 2^{-k} {\hat T}^{k-1}.$$
Докажем, что оператор $\T$  не имеет собственных векторов, кроме констант.  Пусть $\T f=\lambda f.$ 
$\T^nf=\lambda^n f.$
Получим $|P(\lambda)f|=\alpha f$ для некоторого $\lambda$, $|\lambda |=1$. Из
$$|P(\lambda)|= \left|\sum_{k=1}^{\infty} 2^{-k} \lambda^{k-1}\right|=1$$ видно, что   $\lambda=1$.
Но  $T$ -- эргодическое преобразование, значит,  $f$ является константой.

(ii)  Утверждение о том, что спектральная мера $\sigma_T$ преобразования взаимно сингулярна со  
сверточным квадратом $\sigma_T\ast\sigma_T$   
 эквивалентно тому, что   ограничение оператора  $\T$ на пространство $H=Const^\perp$  не имеет 
ненулевого сплетения с произведением $\T\otimes \T$.  Докажем отсутствие ненулевого сплетения.

Пусть
 $J: H \to H \otimes H$ и выполняется условие сплетения
$$J\T = ( \T \otimes \T)J.$$ 

Пусть $\T^{k_j}\to P=P(\T)$, тогда
$$JP = (P \otimes P)J,\  (P \otimes P)^{-1}J =JP^{-1},$$
$$((I-a\T) \otimes (I-a\T))J =(1-a)J(I-a\T),$$
$$[(I\otimes I) +(\T\otimes \T) - (\T\otimes I) - (I\otimes \T)]J=0.$$
 Из равенства
$$\sum_{i,j=0}^{n-1} (\T^i\otimes \T^j)[(I\otimes I) +(\T\otimes \T) - (\T\otimes I) - (I\otimes \T)]J=0$$
получим
$$(I\otimes I)J +(\T^n\otimes \T^n)J - (\T^n\otimes I)J - (I\otimes \T^n)J=0.$$
Свойство слабого перемешивания преобразования $T$ эквивалентно сходимости
$\T^{n_i}\to_w 0$ ( на пространстве $H$) для некоторой последовательности $n_i\to\infty$.
Из сказанного вытекает  $(I\otimes I)J =0,$ $ J=0.$
Что и требовалось.

Теперь приведем  доказательство более сильного утверждения: 
$\hat T\otimes \hat T$ имеет однородный спектр кратности 2. 
Из него также следует упомянутая дизъюнктность,
так как в случае общей компоненты у мер $\sigma_T$  и 
$\sigma_T\ast\sigma_T$ максимальная спектральная кратность
$\hat T\otimes \hat T$ будет  не меньше 4.

Пусть $f$ -- циклический вектор для оператора $\T$, $a=\frac 1 2$.  Рассмотрим векторы 
$$W_n=(I-a\T)^nf\otimes (I-a\T)^nf, \ n\in\N.$$
Обозначим через 
 $C_n$  циклическое пространство с циклическим вектором $W_n$.
Слабые пределы степеней преобразования образуют полугруппу, следовательно, содержат операторы $P^m$.
Таким образом,   для всех $m$  выполнено
$$P^m(I-a\T)^nf\otimes P^m(I-a\T)^nf\in C_n.$$
Отсюда следует, что  для  $k=0,1,2,\dots,n,$ 
$$(I-a\T)^kf\otimes P^m(I-a\T)^kf\in C_n,$$
но это  приводит к 
$$\T^kf\otimes f     +   \T^kf\otimes f=: V_k \in C_n.$$
Так как векторы $V_k$ и их сдвиги  порождают все  пространство $L_2\odot L_2$ симметричных функций $F$ 
 (речь идет о симметрии $F(x,y)=F(y,x)$), мы установили, что циклические пространства $C_n$ 
 приближают   все пространство $L_2\odot L_2$. Значит, оно циклическое.
Оператор $\T\otimes \T$, ограниченный на  $L_2\odot L_2$, имеет простой спектр. Произведение  $\T\otimes \T$ на 
 $L_2\otimes L_2$ имеет  однородный спектр кратности 2, так как  $\T\otimes \T$ есть  прямая  сумма его ограничений
на пространство симметричных функций и пространство антисимметричных функций. А эти ограничения изоморфны друг другу, 
что легко следует из спектрального представления оператора $\T\otimes \T$.

Гипотеза: для  слабо перемешивающих ограниченных конструкций  
выполнено свойство минимальности спектра (MS): степени  $\T^{\odot n}$  имеют простой спектр для всех $n$.

(iii)  
Воспользуемся связью между кратным перемешиванием и самоприсоединениями.
Переходя к подпоследовательности по $i$, обозначая ее также через $i$, получим
 для любых $A,B,C\in\B$
$$\mu(A\cap T^{m(i)}B\cap T^{n(i)}C)\ \to \nu(A\times B\times C),$$
где $\nu$  -- самоприсоединение порядка 3, с проекциями на грани в кубе $X\times X\times X$,  равными $\mu\times\mu$.
  Такая  мера $\nu$  однозначно связана с оператором $J$ равенством
$$(J\chi_A , \chi_B\otimes \chi_C)_{L_2(\mu\times\mu)} =\nu(A\times B\times C),$$
при этом выполняется условие сплетения 
$$J H \subset H \otimes H, \ \  J\T = ( \T \otimes \T)J.$$   
В силу пункта (ii) получим, что
  $J H =\{0\}$,  следовательно, 
$$J L_2 =\{Const\}, \ \ \nu=\mu\times\mu\times\mu,$$
 откуда получаем требуемую ходимость
$$\mu(A\cap T^{m(i)}B\cap T^{n(i)}C)\ \to \mu(A)\mu(B)\mu(C).$$

В \cite{2012}  было дано обобщение результата \cite{JRS}.

\vspace{2mm}
\bf Теорема 4.1. \it  Все нежесткие вполне эргодические  
ограниченные конструкции обладают минимальными самоприсоединениями. \rm

Доказательство этой теоремы основано на эффекте запаздывания.  Когда множество проходит через надстройки,
некоторые его части начинают отставать от других.  Это приводит к тривиализации самоприсоединения. 
Для эргодического самоприсоединения $\nu\neq \Delta_{T^n}$, найдутся $s\neq 0$, $n(j)\to \infty$ и последовательность
множеств $C_j$,  $\mu(C_j)>c>0$
такие, что одновременно 
$$\nu(A\times B)=\lim_j  \mu(T^{n(j)}A\cap B\ | C_j)$$
и 
$$\nu(A\times B)=\lim_j  \mu(T^{n(j)-s}A\cap B\ |C_j).$$
Тогда 
$$\nu(A\times B)=\nu(T^{-s}A\times B), $$
что  в  силу эргодичности $T$  влечет за собой 
$$\nu=\mu\times\mu.$$
Похожий аргумент использовался  в работе \cite{JRS}.
  Алгоритм нахождения множеств $C_j$ описан в \cite{2012}.  

(iv) Покажем, как доказывается MSJ для конструкции $\bar s_j(0,1)$ (рассуждения в общем  случае похожи).

Перепишем (2.1) в виде
$$\nu(A\times B) =\lim_j \mu(T^{k_j}A\cap B\cap Y_j^{k_j}) + \dots. $$
\\
\it Случай 1. \rm Пусть выполнено  $$\frac {k_j} {h_j} \to 0, \ k_j\to\infty, \ \mu( Y_j^{k_j})\to 1,$$
 имеем
$$\nu(A\times B)=\lim_j  \mu(T^{k_j}A\cap B).$$   
Найдем $i(j)$ такие, что
$$h_{i(j)}\leq k_j <  h_{i(j)+1}=2h_{i(j)} +1.$$
Рассмотрим башню этапа ${i(j)}$  и  некоторый этаж $F$  во второй колонне этой башни.  
Тогда образ $T^{k_j}F$  будет состоять из частей, располагающихся на разных этажах в башне:
половина образа на одном этаже, четверть этажом ниже и т.д. Это легко понять, если представить то, как 
этаж движется под действием преобразования $T$, пройдя через надстройки над второй колонной.
Таким образом, на части фазового пространства $C_j$ наблюдается запаздывание на единицу времени и это приводит к

 $$\nu(A\times B)=\lim_j  \mu(T^{k_j}A\cap B)=\lim_j  \mu(T^{k_j}A\cap B | C_j)= \nu(T^{-1}A\times B),$$
 откуда получаем $\nu=\mu\times\mu.$
\\
\it Случай 2. \rm Пусть $\frac {k_j} {h_j} > a>0$   и для некоторого маленького числа $d>0 $ и 
бесконечного множества индексов $j$ выполняется 
$$(1+4d)h_{i(j)}\leq 4 k_j <  (2-4d)h_{i(j)}.$$ 
  Определим $C_j$ как объединение первых $dh_{i(j)}$ этажей.   Меры множеств больше  положительной константы,
считаем, выбирая подпоследовательность, что они сходятся к $c>0$.
Заметим, что 
$$T^{k_j}C_j,\, T^{2k_j}C_j\subset Y_j^{k_j},$$
но  $\mu(T^{k_j}A\cap B \ | Y_j^{k_j})\to \nu(A\times B)$,  поэтому 
$$\mu\left(T^{k_j}(A\cap C_j) \cap B )\ \to с\nu(A\times B\right) \eqno (4.1)$$
и
$$\mu\left(T^{k_j}(A\cap T^{h_{i(j)-2}}C_j) \cap B \right)\ \to с\nu(A\times B). \eqno (4.2)$$
Но одновременно с последней сходимостью будет выполняться 
$$\mu\left(T^{k_j}(A\cap T^{h_{i(j)-2}}C_j) \cap B )\ \to с\nu(T^{-1}A\times B\right). \eqno (4.3)$$
Нам надо докать (4.3). Пусть  $A_j,B_j$ обозначют пересечения  $A$ и $B$, соответственно, 
и  первых колонн этапов $i(j)-2$ и $i(j)-1$. Тогда
$$A=A_j\uu T^{h_{i(j)-2}}A_j\uu\dots,$$
$$B=B_j\uu T^{h_{i(j)-2}}B_j\uu T^{2h_{i(j)-2}+1}B_j\uu\dots.$$
Имеем
$$\mu\left(T^{k_j}(A\cap T^{h_{i(j)-2}}C_j) \cap B \right)=
\mu\left(T^{k_j}(T^{h_{i(j)-2}}A_j\cap T^{h_{i(j)-2}}C_j) \cap T^{2h_{i(j)-2}+1}B_j\right)=$$

$$=\mu\left(T^{k_j}(T^{-1}A_j\cap T^{-1}C_j) \cap T^{h_{i(j)-2}}B_j\right)
=\mu\left(T^{k_j}(T^{-1}A\cap T^{-1}C_j) \cap B\right).$$
Подставив в (4.1)   множество $T^{-1}A$ вместо $A$, в виду очевидной сходимости $\mu(C_j\Delta T^{-1}C_j)\to 0$
устанавливаем (4.3):
$$\mu\left(T^{k_j}(T^{-1}A\cap T^{-1}C_j) \cap B)\to c\nu(T^{-1}A\times B\right).$$ 
С учетом  (4.2)   получаем
$$ \nu(A\times B)= \nu(T^{-1}A\times B),\ \nu=\mu\times\mu.$$

Если случай 2 не реализуется, то 
$h_{i(j)-1}\approx  k_j$   (или $h_{i(j)-2}\approx  k_j$). 
Замечаем, что теперь мы  оказываемся в условиях случая 1 с новой последовательностью  $k_j'=k_j-h_{i(j)-1}$:
$$\mu(T^{k_j'}A\cap  B)\to \nu(A\times B).$$

Таким образом,   показано, что единственными эргодическими присоединениями могут быть 
только $\mu\times\mu$ и $\Delta_{T^k}$
(когда $k_j=k$). Свойство MSJ установлено.

\bf Дизъюнктность степеней. \rm Из  MSJ  вытекает  дизъюнктность положительных 
степеней преобразования: $\mu\times\mu$  является   единственным  
присоединением степеней  $T^m$ и $T^n$  при  $m,n>0, m\neq n.$
Эквивалентное утверждение: если $\T^qJ=J\T^p$, $q\neq p$,  где $J$ -- марковский оператор, то  $J=\Theta$.
В \cite{2012}  в связи с работой \cite{B1} доказано следующее утверждение. 

\vspace{3mm}  
\bf Теорема 4.2.  \it  Для всех слабо перемешивающих   
ограниченных конструкций $T$  при $0<m<n$   степени $T^m$ и $T^n$ дизъюнктны. \rm

Если слабо перемешивающая конструкция не обладает жесткостью, то она обладает MSJ.
Для жесткой слабо перемешивающей ограниченной конструкции $T$ найдется   сколь угодно малое число   $\eps>0$ 
такое,  что 
$$\T^{qh_j}\to_w Q(\T)=q\eps I+(1-q\eps)R(\T)$$   и $$\T^{ph_j}\to_w P(\T)=p\eps I+(1-p\eps)R(\T),$$ 
где ряд $R(\T)\neq I$ ( см. \cite{2012}).
Дизъюнктность  $T^q$  и  $ T^p $ при $p>q$ очевидна, так как мы получаем $R(\T)J = J$, что влечет за собой $J=\Theta$.

Другой удобный инструмент для доказательства дизъюнктности степеней преобразований дает лемма из \cite{mmo13}.

\vspace{3mm}
\bf  Лемма. \it Пусть $S$, $T$ —  вполне эргодические  пребразования такие, что
$$\hat S^qJ = J\T^p, $$
где $q, p$ взаимно просты, а $J\neq \Theta$ --
неразложимый марковский оператор. Если для полиномов $Q, P$ выполнено
$$Q(\hat S)J = JP(\T),$$ 
то найдется полином $R$ такой, что 
$$Q(\hat S) = R(\hat S^{q}), \ \ P(\T) = R(\T^{p}).$$ \rm

Результат  о спектральной дизъюнктности степеней ограниченных и некоторых неограниченных конструкций 
был получен  в работе \cite{ALR}.


\section{Полуограниченные конструкции}
Рассмотрим конструкции с надстройками  вида $ \bar s_j=(0,s_j,0)$. Любопытно, что в случае
 $s_j=[ j^a]$ при $a=0$ централизатор
тривиален (модифицированное преобразование Чакона), при $0<a<1$, как мы покажем,   централизатор континуален, но
   при $a\geq 1$ он снова тривиален.

Скажем,  что  $s_j$  медленно растет,  если она монотонная,    каждое значение
 $r$ она принимает $N(r)$ раз, причем 
$N(r)\to\infty$ при $r\to\infty$. Примеры: $s_j\sim \sqrt j$,  $s_j\sim \ln j$.

\vspace{3mm}
\bf Теорема 5.1. \it  Пусть $T$ -- конструкция  с надстройками вида $\bar s_j=(0,s_j,0)$, где $s_j$  --
медленно растущая последовательность. Тогда все ряды вида 
$$P(\hat T)=\sum_{k\in \Z}a_k {\hat T}^{k},\ a_k\geq 0,\ \sum_{k\in \Z}a_k =1, $$
лежат в слабом замыкании степеней оператора  $\T$.

\bf Следствие. \it Преобразование $T$ является жестким,
 т.е. $\hat T^{m_i}\to I$ для  некоторой последовательности $m_i\to\infty$.
 Cпектры  симметрических  тензорных степеней $\hat T^{\odot n}$  являются простыми.\rm
 
Доказательство. Нам  будет удобнее писать $h(j)$ вместо $h_j$. Обозначим через $j[k]$
минимальное число $i$ такое, что $k=s_i$. С учетом того, что $s_{j[k]+m}=k$ для $m=1,2,\dots N(k)-1$,
получим для больших значений $k$
$$\hat T^{h(j[k])}\approx_w  \frac 1 2  I+  \frac 1 2 \T^{-k}.$$
Теперь  подставим  $h_j$ (в других случаях  $-h_j$) вместо  $k$.  Получим
$$\hat T^{h(j[j[k]])}\approx_w  \frac 1 2  I+  \frac 1 4 I +  \frac 1 4 \T^k, $$
$$\hat T^{h(j[j[j[k]]])} \approx_w  I+  \frac 1 4 I +  \frac 1 8 I +  \frac 1 8  \T^{-k} $$
и так далее.  
Находим последовательность $n_i$, для которой выполнено 
$$\hat T^{n_i}\to_w \frac 1 2  I+  \frac 1 4 I +  \frac 1 8 I + \dots = I. $$
Так как 
$$\hat T^{n_i+m}\to_w   \T^m,  $$
получаем 
$$\hat T^{-h[n_i]}\approx_w  \frac 1 2  I+  \frac 1 2 \T^{n_i}\approx_w  \frac 1 2  I+  \frac 1 2 \T^m.$$
Если слабоe замыканиe степеней эргодического преобразования 
содержит пределы вида $\frac 1 2  (I+  \T^m)$, то замыкание содержит все допустимые полиномы
( доказательство  см. в \cite{07}).
В силу сказанного спектры  симметрических  тензорных степеней 
$\hat T^{\odot n}$  являются простыми.  Преобразование $T$ является жестким,
 ( $\hat T^{n_i}\to I$ для  некоторой последовательности $n_i\to\infty$), это принципиально   отличает
такую конструкцию от преобразования Чакона с $s_j=1$, обладающего свойством минимальных самоприсоединений.

Вопрос: будет ли  $(0,[\sqrt j],0)$-конструкция обладать свойством простоты в смысле теории
присоединений? В этом случае ее централизатор континуальный, 
но при  $s_j=j$   централизатор тривиален. Покажем это.
 Если 
$$\hat T^{n_i}\to I, \ \ h_{j(i)}\leq n_i<h_{j(i)+1},$$
то из-за  надстроек части пространства смещены на величины надстроек, что приводит к  
$$\hat T^{n_i-j(i)}\to I,\  \hat T^{n_i-j(i)-1}\to I,$$
сходимость к $I$ выполняется также для разности степеней, поэтому получим
но тогда $\T=I$. 

А в случае   $s_j=j^2$ получим 
$$\hat T^{n_i-(j(i)+1)^2}\to I, \hspace{2mm}  \hat T^{n_i-(j(i)+2)^2}\to I,{}\,\,\,\,\, \  \hat T^{n_i-(j(i)+3)^2}\to I,$$
$$\hat T^{-(j(i)+1)^2 +(j(i)+2)^2}\to I,{} \ \  \hat T^{-(j(i)+2)^2+(j(i)+3)^2}\to I,$$
$$\hat T^{2j(i)+3}\to I, \ \ \hat T^{2j(i)+5}\to I,$$
$$\T^2=I.$$

\bf  Пример конструкции с фактором. \rm  Рассмотрим преобразование, заданное  последовательностью
$$\bar s_j=(0,2^j,0,0).$$
Имеем
$$h_{j+1}=4h_j+2^j,$$

$$\frac {2^{j+1}} {h_{j+1}} < \frac {2^j} {2h_j}.$$

Следовательно,
$$\sum \frac {2^j} {h_j}<\infty.$$
Это  означает, что рассматриваемое преобразование действует на пространстве с конечной мерой.
Он содержит фактор, изоморфный одометру с двоично-рациональным спектром (с параметрами $r_j=2$).
Действительно, начиная с этапа $j$,  все параметры и высота  $h_j$  будут кратны $2^j$. 
Следовательно, в спектр содержит группу  $\{e^{\frac{\pi i k} {2^j}}\}$.

Читатель может рассмотреть различные последовательности   $s_j$,  например, $j, j^n$,$2^n$  и  $p(j)$, где $p(j)$ есть
$j$-тое простое число.   Какие спектральные свойства у соответствующих конструкций, 
есть ли у них факторы, является ли 
централизатор тривиальным?  На некоторые из них ответы были даны  в \cite{2013c}.  

 \bf  Стохастическое  преобразование Чакона. \rm  Пусть  $r_j=j$, рассмотрим всевозможные последовательности 
   $(s_j(1), s_j(2),\dots, s_j(j-1),s_j(j))$, где высоты $s_j(i)$ надстроек принимают с одинаковой вероятностью 
независимо  значения 0 и 1. В результате мы имеем ансамбль конструкций ранга один. 
В работе  \cite{07} была высказана гипотеза:
 
\it  для почти всех  стохастических конструкций $T$  слабое замыкание их степеней 
 состоит из $\Theta$ и операторов вида $T^sP^m$, где $P=\frac 1 2 (I+\T)$. \rm

Полиномиальная структура слабого замыкания   
  с учетом известных фактов о  преобразованиях ранга один  влечет за собой 
свойство минимальных саоприсоединений.  

В связи с упомянутой гипотезой  возникла  другая задача,
представляющая  самостоятельный интерес.

Пусть  $f: \Z_r \to \{0,1\}$. Положим  
$$P(f, m, s)=\left|\left\{i\in \Z_r \ : \ \sum_{w=1}^m f (z + w)= s \right\}\right|,$$
где  $z$ и $w$ складываются по модулю $r$, и определим
$$D(f,m)=\sum_{s=1}^m \left| P (f, m, s) - P (f, m, s-1) \right|.$$

\bf Cтатистическая лемма. \it Для любого $ \eps > 0$ найдется  натуральное $L$ такое, что для 
всех больших  значений $r$ типичная функция $f: \Z_r \to \{0,1\}$  
   удовлетворяет  условию $$ D(f,m)< \eps r, \ \  L< m <r - L. $$\rm

Из леммы втекает  гипотеза (схема доказательства изложена в  \cite{2012}). М.Е. Липатов сообщил
автору факт из теории случайных блужданий, при помощи которого можно доказать лемму. 
Этому вопросу планируется посвятить  отдельную заметку.
  
Отметим, что преобразование с указанной полугруппой слабых пределов  можно получить модификацией стохастического
преобразования Чакона, не используя леммы. Грубо говоря, новые элементы конструкции обеспечат те свойства, 
которые вытекают из леммы. Новые элементы -- это редкие высокие надстройки орнстейновского толка, 
которые не  оказывают влияния  на тех участках времени,
где формируются пределы  $T^sP^m$, но обеспечивают близость к $\Theta$  на остальных временных участках. 
При этом на некоторой последовательности этапов используются только обычные орнстейновские надстройки.


\section{ Стохастические конструкции} 
Д.Орнстейн доказал свойство перемешивания для большинства конструкций
определенного им  статистического ансамбля \cite{O}.  Напомним его определение.
Фиксируем  последовательность $r_j\to\infty$  и  $H_j\to \infty$,  $H_j<<r_j$,  
рассмотрим  всевозможные
 последовательности вида $a_j(i)\in \{0, 1, \dots, H_j-1\}$.
Пространство $$  \{0, 1, \dots, H_1-1\}\times \{0, 1, \dots, H_2-1\}\times\dots$$
оснащается   естественной вероятностной мерой (произведение равномерных распределений).
Кострукции задаются параметрами $$ s_j(i)=b_j+a_j(i)- a_j(i+1),$$
 далее считаем, что  $b_j=H_j$.

 Свойства конструкции $T$  зависят от свойств последовательности $a_j(i)$.
Какие достаточные условия гарантируют свойство перемешивания $T$ ?
Типичная последовательность $a_j(i)$ при фиксированном $j$ принимает значение $s\in \{0,1,\dots, H_j-1\}$
с частотой, близкой к  $\frac 1 {H_j}$,   а значения разностей $a_j(i)- a_j(i+1)$ распределены треугольно:
доля $i$, для которых выполнено  $a_j(i)- a_j(i+p)=s$ при  $s\in \{-H_j+1, -H_j+2,\dots,H_j-2, H_j-1\}$,
близка к $\frac {H_j - |s|}{H^2_j}$. Такая ситуация влечет за собой свойство перемешивания конструкции
даже в том случае, когда распределению частот асимптотически разрешено отличаться от  от треугольного в любое,
 но фиксированное число раз.

Обозначим 
$$S_j(i,p)=s_j(i)+s_j(i+1)+\dots+s_j(i+p-1),$$
$$Q_{j,p}=\frac 1 {r_j -p}\sum_{i= 0}^{r_{j}-p -1}\T^{-S_j(i,p)}$$
 Укажем  условия, обеспечивающие  свойство перемешивания.  

\vspace{3mm}
\bf Теорема 6.1. \it  Если для любого $\eps>0$  для всех достаточно больших $j$ для всех $p$, 
$0<p< (1-\eps)r_j$  выполняется $ Q_{j,p}\approx_s \Theta, $ 
то конструкция $T$  обладает перемешиванием. \rm
    
Схема доказательства.   Пусть  $h_j\leq m <h_{j+1}$, тогда 

$$m= ph_j + q + \sum_{i=1}^{p} s_j(i), \  0\leq q<h_j.$$
Специфика конструкций ранга один позволяет оценивать величину $\mu(T^{m}A\cap B)$ методом трех областей.
Область $D_0$ является объединением колонн (этапа $j$)  с номерами от   1 до $p$,
 область $D_1$ является объединением верхних частей  оставшихся  колонн (этажи  с номерами $n$ при $q<n<h_j$).
$D_2$  -- оставшееся множество,  соответствующее объединению нижних частей колонн (этажи  с номерами $n$ при $0<n<q$). 

Пусть $A, B$ состоят из этажей башни $\xi_j$.  Тогда имеем
$$\mu(T^{m}A\cap B)\approx \frac 1 {r_{j -1}}\sum_{i= 1}^{r_{j+1} -1}
\mu(T^{m_j-h_{j+1}}T^{-s_{j+1}(i)}A\cap B\cap D_1)+$$ 
$$ 
 +\frac 1 {r_j -q}\sum_{i= 1}^{r_j -q}
\mu(T^{k}T^{-S_j(i,q)}A\cap B\cap D_2)\ +
 \frac 1 {r_j -q-1}\sum_{i= 1}^{r_j -q-1}
\mu(T^{-h_j + k}T^{-S_j(i,q+1)}A\cap B\cap D_3).
$$

Если выполняется достаточное условие, то 
$$(f,T^mg)\approx \int_D f\T^k Q_{j,1}g\,d\mu \ +\ \int_{D_1} f\T^{k_1} Q_{j,p}g\,d\mu+ $$
$$+\int_{D_2}f \T^{k_2} Q_{j,p+1}g\,d\mu \approx (f,\Theta g)=\int f g\,d\mu \int f g\,d\mu.
$$
Это означает, что наша конструкция перемешивает.

Осталось проверить, что достаточное условие выполнено.

В нашем случае мы имеем
$$S_j(i,p)= pH_j + a_j(i)-a_j(i+p),$$
поэтому
$$\T^{- pH_j}Q_{j,p}=\frac 1 {r_j -p} \sum_{i= 0}^{r_{j}-p -1}\T^{ a_j(i)-a_j(i+p)}.$$

Но
$$
\left\| \T^{ pH_j}Q_{j,p} -  P_j^\ast P_j \right\|\approx 0,
$$
где
$$P_j=\frac 1 {H_j}\sum_{i=0}^{H_j-1}T^i.$$

Хорошо известно, что для эргодических $T$ выполнено 
$$
P_j\approx_s \Theta, \ \   P_j^\ast  P_j\approx_s \Theta. 
$$
Получили, что 
$$
Q_{j,p}\approx_s \T^{ pH_j}\Theta =\Theta. 
$$

Поясним, в каком  смысле распределение значений $S_j(i,p)$ может  отличаться  
от  равномерного или  треугольного, но при этом обеспечивает свойство перемешивания конструкции.

\vspace{3mm}
\bf Лемма 6.2. \it  Пусть $R_j$  -- некоторая последовательность марковских операторов, 
коммутирующих со слабо перемешивающей
конструкцией $T$, $R_j\to_s \Theta.$
Если, для некоторого $a>0$  при  $0<p< (1-\eps)r_j$
для всех достаточно больших $j$ выполнено
$R_j\geq a Q_{j,p},\ \ R_j\to_s \Theta,$ то конструкция $T$ перемешивает.\rm

Доказательство.  $$R_j^\ast R_j\to_s \Theta,  $$ 
$$ a^2Q_{j,p}^\ast Q_{j,p}+\dots \to_s  \Theta. $$
Так как оператор  $\Theta$ неразложим в марковском централизаторе оператора $\T$ ($\mu\times\mu$ эргодична),
а $Q_{j,p}$ лежат в централизаторе, получаем
$$ Q_{j,p}^\ast Q_{j,p}\to_w  \Theta, \ \  Q_{j,p}\to_s  \Theta.$$

\section{Типичность ранга один, нестандартное перемешивание, нетипичные преобразования}
В \cite{Ba} доказана типичность ранга один в пространстве Mix  на основе результата А.И. Баштанова о
близости любого класса сопряженности  к бернуллиевскому преобразованию и теоремы С.В.Тихонова о плотности
 бернуллиевских преобразований в Mix.   Типичность ранга один означает следующее:
множество перемешивающих преобразований ранга один имеет тип $G_\delta$ и оно плотно в Mix.

Наш подход другой: мы 
 непосредственно  покажем, что  класс сопряженности 
орнстейновской   конструкции плотен в пространстве Mix.
Напомним, что преобразования $S$ и  $T$ близки в Mix, если для всех $k$  степени $\hat S^k$ и  $\T^k$  
близки в слабой операторной топологии.  Заметим, что для случая пространств 
с бесконечной мерой годится то же самое определение.

Пусть $S$ -- перемешивающее преобразование, а $T'$  заранее фиксированная орстейновская конструкция.
Надо найти такое сопряжение $R$, чтобы не только  $T=R^{-1}T'R$ было близко к $S$ в слабой тополоии, 
но и то же самое 
было верным для всех их степеней.  Это означает, что, когда  
степени $S^i$  станут  близки к оператору $\Theta$, мы можем забыть про них, переключив внимание на то,
 чтобы степени $T^i$ также были близки к $\Theta$.   Это и будет означать, что преобразование $S$ и 
сопряжение $T$ близки  в Mix.

Выбираем башню Рохлина-Халмоша преобразования  $S$  высоты $N+a$ с малым остатком.
 и башню $\xi'$  преобразования $T'$  ранга один   высоты  $h=h_j$ для очень далекого этапа $j$.
    Считаем, что $h=Nq$,    $q>>N>>a>> 1$. 
В башне $\xi'$  отметим этажи с номерами $iN +a(i)$,  $i=0,1,2, \dots, q-1.$
Теперь намотаем башню $\xi'$ на башню  $\xi$ преобразования $S$ так, что отмеченные этажи окажутся на первом этаже
$\xi$, а куски орбит длиной 
$N+ a(i+1)-a(i)$ (между отмеченными этажами) совместились с кусками орбит преобразования $S$ той же длины.

Намотка -- это и есть сопряжение, причем мы  получаем сопряженое преобразование  $T$, совпадающее с $S$ 
 на $N-a$  первых этажах.  Легко реализовать  ситуацию, когда  степени 
$T^i$, $S^i$ мало отличаются для всех значений $0<i<M$, а при $i>M$ выполнено $S^i\approx\Theta$.
Мы достигаем поставленной цели, если  убедимся в том, что $T^i\approx\Theta$ для всех $i>M$. 
Посмотрим, к примеру, что мы имеем при $i\approx N/2$.  На верхней половине $D$  башни $\xi$  степень $T^i$  совпадает
 с  $S^i$ (на маленькое множество
меры, сравнимой с $a/N$ не обращаем внимания). Но   $S^i\approx_w\Theta$, поэтому на этой половине степень $T^i$
близка к $\Theta$.
Это надо так понимать: для заранее фиксированных множеств $A,B$ 
 $$\mu( T^iA\cap B\cap D)\approx \mu(A)\mu(B) \mu(D).$$
А на нижней половине перемешивание  обеспечиваеся  виртуальными надстройками 
$a+ a(i+1)-a(i)$. Тут нужно учесть, что
заранее фиксированное множество равномерно распределено в башне $\xi$, а величина $q$   в соотношении  $h=Nq$
 выбрается  весьма большой, $q>>2^N$. При $i>N$ мы забываем про $S^i$ окончательно,  степени  $T^i$ перемешивают
заранее выбранный фиксированный набор  множеств благодаря виртуальным, а потом и реальным надстройкам.
Такова  вкратце схема доказательства плотности орнстейновского преобразования ранга один в Mix.  Из этого также
можно получить утверждение о типичности ранга один в пространстве Mix.

  Аналогично  доказывается плотность орнстейновского потока в пространстве перемшивающих потоков.

\bf Теорема 7.1. \it В пространстве перемешивающих потоков типичными являются потоки ранга один.  Следовательно, 
типично свойство минимальных самоприсоединений.\rm

В случае пространств с бесконечной мерой доказательство типичности ранга один можно упростить за счет привлечения
идеи сидоновских надстроек.

{ \bf К вопросу об изоморфизме   перемешивающих конструкций.}  В \cite{11dis} дано достаточное  
условие того, что две перемешивающие 
конструкции  $T$, $T'$ в случае вероятностного пространства  не изоморфны.  

\bf Теорема 7.2. \it Если у перемешивающих конструкций ранга один $T$, $T'$  высоты 
$h_j$, $h_j'$ башен сравнимы,  но их отношение не стремится к 1, то $T$, $T'$ неизоморфны. \rm

Аналогичное утверждение справедливо для перемешивающих потоков ранга 1. Сформулируем следствие из \cite{97}.

\bf Теорема 7.3. \it  Перемешивающий поток $T_t$ ранга один не изоморфен потоку $T_{at}$ при $a>1$. \rm

\vspace{3mm}
\bf Нестандартное перемешивание. \rm  В случае, когда отношение высот башен стремится к 1, 
интересно  найти  некоторые  свойства, инварианты,  которые бы отличали 
одну перемешивающую конструкцию от другой.

Пусть $T$ -- перемешивающее преобразование ранга один с соответствующей последовательностью разбиений $\xi_j$.
Пусть множества $A_j$,  $B_j$ состоят из атомов разбиения  $\xi_j$, причем  эти множества асимптотически 
 равномерно распределены в башне.
Это означает что большинство интервалов  из $L=L(\eps)$  последовательных этажей  с точностью до $\eps L$ содержат 
$\mu(A_j)L$ этажей,  лежащих в $A_j$. Такие последовательности $A_j$ назовем правильными. (При работе с 
такими последовательностями можно наложить  условие $\mu(A_j)\to a>0$.)

Преобразование \it хорошо перемешивает, \rm если для  любых правильных $A_j$, $B_j$ выполнено
$$ \sup_{m>h_j}\mu(T^mA_j\cap B_j)-\mu(A_j)\mu(B_j)\ \to 0,  \ j\to \infty.$$
Это определение  навеяно нестандартным анализом: здесь преобразование перемешивает не только стандартные множества
($A_j=A$, $B_j=B$), но и некоторые нестандартные множества. 

Почти все  орнстейновские преобразования  хорошо перемешивают.
А лестничные конструкции нет. Для них  реализуется случай $\mu(A_j)\to \frac 1 2 $ и
$$\mu(T^{2h_j}A_j\cap A_j)-\mu(A_j)^2\ \to \frac 1 4.$$
В качестве таких  $A_j$ можно взять объединение нечетных этажей в башне разбиения $\xi_j$. 

{ \bf Нетипичные свойства.}  Следующий вопрос имеет смысл  как  для пространства Aut, 
так и для пространства  Mix:
\it  Пусть $K$ -- компактное множество преобразований. Верно ли что найдется плотное $G_\delta$-множество  $Y$, 
не пересекающееся с $K^{Aut}$?\rm

Через $K^{Aut}$ мы обозначаем объединение 
всех классов сопряжености в группе ${Aut}$, имеющих представителя в $K\subset Aut$.
Фиксируем метрику  $dist$, задающую слабую операторную топологию. 

\vspace{3mm}
\bf Теорема 7.4.  \it Пусть $K\subset Aut$ -- компакт и для некоторого $r>0$ для всех 
преобразований $T\in K$ и любого натурального числа $n$ найдется $ m>n$ такое, что  $dist(T^m, \Theta)>r$.
Тогда найдется всюду плотное $G_\delta$-множество  $Y$, не пересекающееся с $K^{Aut}$. \rm
  
Доказательство. Пусть  $m(T,j)$ -- минимальное число среди тех $m>n$, для которых  $dist(T^m, \Theta)>r$.
Так как $K$ -- компакт, $m(T,j)$ -- ограниченная функция на $K$. 
Обозначим через $M(j)$  максимум значений $ m(T,j)$.  Рассмотрим множества $F_j=\{j,j+1,\dots, M(j)\}$.
В \cite{U06} доказано, что существует  всюду плотное $G_\delta$-множество  $Y$ такое, что для каждого $S\in Y$
   найдется перемешивающая подпоследовательность $F_{j(k)}$, $j(k)\to\infty$, т.е. при  $m(k)\in F_{j(k)}$
$$ dist(S^{m(k)}, \Theta)\to 0, k\to\infty. $$
Такое преобразование $S$ и любое, сопряженное ему,  не принадлежит $K$, так как иначе для всех  $k$ 
должно выполняться
$$\exists m\in F_{j(k)} \ \ dist(S^m, \Theta)>r.$$
Теорема  доказана.

Можно привести разнообразные  примеры  компактов, удовлетворяющих условию теоремы. 
Простейший случай:  конечный набор неперемешивающих преобразований. Хорошо известные примеры  --
 перекладывания $k$ отрезков. Как следствие, мы  получаем результат  \cite{CD} о том, что перекладывания
конечного числа отрезков не  являются типичными преобразованиями. Наше доказательство 
короче, так как мы не заботились об оценке чисел $M(j)$.  Условие  $dist(T^m, \Theta)>r$  для перекладываний
выполняется в силу хорошо известного свойства частичной жесткости. 

Вопрос: может ли перекладывание конечного числа прямоугольников  обладать свойством перемешивания?


\section{ Алгебраические параметры вместо случайных}

Пусть  $r_j$  -- некоторая последовательность простых чисел,  $r_j\to\infty.$ 
Выберем $q_j$ -- образующий элемент мультипликативной группы  поля $\F_{r_j}$, которое отождествляется 
с набором натуральных чисел $0,1,\dots, r_j-1$.    Положим  
 $$ s_j(i)=r_j+\{q_j^i\} - \{q_j^{i+1}\},$$
где $\{q\}$ -- означает натуральное число (остаток),  отвечающее элементу поля $q$.
Теперь мы рассмотрим конструкцию $T$   ранга один с такими параметрами. Нам надо установить свойство перемешивания.
Положим
$$S_j(i,p)=s_j(i)+s_j(i+1)+\dots+s_j(i+p-1),$$
$$Q_{j,p}=\frac 1 {r_j -p}\sum_{i= 0}^{r_{j}-p -1}\T^{-S_j(i,p)}.$$
 
Имеем
$$S_j(i,p)= pr_j + \{q_j^i\} - \{q_j^{i+p}\},$$

$$\T^{- pr_j}Q_{j,p}=\frac 1 {r_j -p} \sum_{i= 0}^{r_{j}-p -1}\T^{ \{q_j^i\} - \{q_j^{i+p}\}}.$$

Но функция $\{q_j^i\} - \{q_j^{i+p}\}$ как функция от $i$  является инъективной:
$$\{q_j^i\} - \{q_j^{i+p}\}=\{q_j^k\} - \{q_j^{k+p}\},$$
$$q_j^i - q_j^{i+p}=q_j^k - q_j^{k+p},$$
$$q_j^i=q_j^k,   \ \ i=k. $$

Положим
$$R_j=\frac 1 {r_j -p}\sum_{i=pr_j}^{(p+1)r_j -p}T^i.$$
Так как $R_j\approx_s \Theta$, чтобы установить  свойство перемешивания,  осталось проверить, что конструкция обладает свойством  слабого перемешивания 
и применить  лемму 6.2. 

Приведем  пример конструкции  М.С.Лобанова и автора.  Рассмотрим отображение  следа 
 $$ tr: \F_{b^n}\to\F_{b}, \ \ tr(q)=\sum_{s=0}^{n-1}q^{b_j^s}.$$
Положим $$r_j=b_j^{n_j}-1\to\infty, \ \ a_j(i)=tr(q_j^i)=\sum_{s=0}^{n-1}q_j^{b^s}.$$
 $$ s_j(i)=b_j+\{ a_j(i)\} - \{ a_j(i+1)\},$$
Эти параметры также определяют перемешивающую конструкцию $T$.
В отличие от орнстейновского случая стохастические  последовательности заменены на квазислучайные.
Подставим в  лемму $P_j,R_j$, положив   
$$P_j=\frac 1 {r_j}\sum_{i=0}^{b_j-1}T^i,  \ \   R_j=\T^{- pb_j} P_j^\ast P_j.$$
В этом случае распределение весов у сумм $Q_{j,p}$ подчинено треугольному распределению.
И тем самым мы получаем новые эффективные примеры  премешивающих конструкций.
  Интерес может представить  изучение спектров таких конструкций.

\bf Простой спектр тензорных произведений.  \rm Параметры, заданные при помощи функции следа 
$ tr: \F_{2^n}\to\F_{2}$, 
 использовались в $\cite{LR}$ для  построения потока $T_t$ такого, что произведения 
$T_{1}\times T_{t}$ имеют простой спектр  для любых  $t>1$. 

Аналогичная  конструкция удовлетворяет более сильному условию:  произведения  
$$T_{t_1}\times T_{t_2}\times T_{t_3}\times \dots$$ имеют простой спектр  для любых различных $t_i>0$.
Причина в том, что такие произведения 
 обладают всевозможными слабыми пределами вида
$$P_b\otimes P_b \otimes P_b \otimes\dots, \ \ \ 4P_b=\T_{-b}+2I+ \T_b,  \ b\in\R.$$
Наличие этих пределов обеспечивает простоту  спектра  произведений при условии, что  исходный поток  $T_t$ 
 обладает простым спектром.
В конструкции потока $T_t$ используются 
надстройки вида 
 $$   s_j(i) =a_j(i)-a_j(i+1) + b_j, \ \ a_j(i)=|tr(q_j^i)|b_j,$$ 
где последовательность $b_j$-- бесконечное число раз принимают каждое рациональное значение.
Методы \cite{07} позволяют строить перемешивающий поток  с аналогичным спектральным свойством.


 \section {Лестничные  конструкции}
 Лестничная конструкция  задается последовательностью надстроек 
$$s(i)=i,  \ i=1,2,\dots, r_j,\ \
r_j\to\infty.$$
   T. Адамс \cite{Ad} нашел оригинальный способ доказательства свойства перемешивания 
таких преобразований в случае $r_j^2/h_j\to 0$.

Для лестничной конструкции легко установить перемешивание для  последовательности  $\{m_j\}$, когда
$m_j\in[h_j, Ch_j]$ для фиксированного числа  $C>1$.
Например, используя эргодичность степени $T^p$, получим
$$\mu(T^{ph_j}A\cap B)\approx \frac {1} {r_j -p}\sum_{i=0}^{r_j -p-1}
\mu(T^{-pi-k(p,j)}A\cap B)\approx \mu(A)\mu(B).$$

В общем случае ситуация аналогична той, что мы рассматривали в \S 5, но здесь   возникают  усреднения вида
$$\frac 1{N_j }\sum_{i=0}^{N_j-1 }
\mu(T^{-d_ji +k_j}A\cap B\cap D), $$
где $d_j$, вообще говоря,   не ограничены.  

Таким образом, нам надо  показать, что  
$$P_j =\frac 1{N_j }\sum_{i=0}^{N_j -1}
\T^{-d_ji}\approx_w \Theta.$$
Фиксируем  большое $L$.  Если 
$d_j\in [h_{p(j)}, 2h_{p(j)}],$  то выполнено 
$$id_j\in [h_{p(j)}, 2Lh_{p(j)}].\eqno (9.1)$$
Отсюда вытекает, что 
$$ \T^{d_j}, \T^{2d_j},\dots ,\T^{(L-1)d_j}\approx_w\Theta,\eqno (9.2) $$
обозначив 
$$Q_L =\frac 1 {L }\sum_{i=0}^{L -1} \T^{d_ji},$$
имеем в силу (9.2)
$$Q^\ast_LQ_L\approx_w\Theta,$$
что равносильно 
$$ Q_L\approx_s\Theta.$$
Теперь получаем
$$P_j\approx_s P_jQ_L\approx_s P_j\Theta\approx_s\Theta.$$

Но (9.1) может не выполняться.  Тогда  Адамс 
 находит такое такое число $a$, что  
$$aLd_j <<N_j,  \ \ ad_j\in [h_{p(j)}, 2h_{p(j)}].$$
Это приводит к желаемой аппроксимации
$$A_L =\frac 1 {L }\sum_{i=0}^{L -1}
T^{iad_j} \approx_s\Theta,$$
$$P_j\approx_s P_jA_L\approx_s\Theta.$$

Особенностью метода Адамса является то, что перемешивание на этапе $j$ использует информацию, полученную
на предыдущих этапах.  В доказательстве использовалось ограничение  $r_j^2/h_j\to 0$
на рост последовательности $r_j$.  Это условие позволяет осуществлять контроль перемешивания, используя только 
три усреднения,
отвечающие трем областям фазового пространства. Мы их описывали в \S 5.  В случае   $r_j^2/h_j\to \infty$ 
возникает неограниченное количество областей $D,D_1,\dots ,D_k,\dots, $  
каждой из которых соответствует  оператор усреднения
вида  
$$Q_{L,k} =\frac 1 {L }\sum_{i=0}^{L -1}
\T^{i(d_j+k)}.$$ 
Так как   $T$ --  слабо перемешивающее преобразование, для большинства номеров $k$  выполнено
$$ \T^{d_j+k} \T^{2(d_j+k)}\dots \T^{(L-1)(d_j+k)}\approx_w\Theta, $$
следовательно, для большинства номеров $k$  выполнено $Q_{L,k}\approx_s\Theta. $
Это показывает, что  на большинстве областей  наблюдается эффект  перемешивания, 
следовательно,  конструкция обладает свойством перемешивания.
Подробно этот  подход  изложен в \cite{stair}. Над проблематикой перемешивания лестничных 
конструкций также работали авторы \cite{CS}, 
 получившие ряд обобщений результата Адамса.

Последовательность  $r_j$ может расти так  быстро, что фазовое пространство
лестничной конструкции будет иметь  бесконечную меру.  В следующем утверждении
речь идет как о   пространствах с конечной,  так и пространствах  с бесконечной мерой.

\bf Теорема 9.  \it Если выполнено условие $\frac {r_j} {h_j}\to 0$, лестничная конструкция обладает перемешиванием.
   \rm

Отметим, что случай $r_j\sim h_j$ представляет собой трудную задачу, сходимость  $\T^i\to_w 0$ не доказана,
а упомянутые методы не работают.


\section{Явная перемешивающая конструкция с двукратным  спектром }
 Интерес к лестничным конструкциям и их модификациям был связан с задачей  
об однородном спектре преобразования 
в классе  Mix.   В \cite{1999}
 был указан класс преобразований и доказывалось, что в этом классе  найдутся нужные
конструкции $T$, у которых $T\times T$ имеет двукратный спектр. Какие именно, оставалось неясным.

Положим $\bar s_j=(1,2,\dots, r_j-2,r_j-1,0)$ для всех этапов, начиная с некоторого.   
Будем считать, что   последовательность $r_j$ начиная с некоторого момента  
$2^{4r}$ раз подряд  принимает значение $r$, потом 
$2^{4(r+1)}$ раз подряд $r_j$ принимает значение $r+1$ и т.д. 
Назовем такой рост последовательности $r_j$  медленным. 
Наша цель -- доказать следующее утверждение.

\vspace{3mm}
\bf Теорема 10.1. \it Если последовательность  $r_j$ конструкции $T$ имеет медленный рост, то   
 $\T\otimes\T$ имеет  однородный двукратный спектр.   \rm

Фиксируем индикатор $f$  некоторого этажа.  Известно, что $f$ является циклическим вектором.
Докажем, что для нашей конструкции $T$  для всех $s>0$  векторы $\T^sf\otimes f + f\otimes \T^sf$
принадлежат циклическому пространству $C_{f\otimes f}$  оператора $\tt=\T\otimes \T$.
Это означает, что симметрическая степень $\T\odot\T$ имеет простой  спектр, откуда следует кратность 
2 для $\T\otimes\T$.

Предположим, что  для некоторой последовательности $M(r)\to\infty$, $r\to\infty$, для всех  $m,n$, 
 $1 \leq m<n\leq M(r)$, выполнены неравенства 
$$ \left|\left( \tt^{p(m,r)}F, \Q_rF \right) - \left( \Q_rF, \Q_rF \right)\right|<\eps_r,$$
$$\left| \left( \tt^{p(m,r)}F, \tt^{p(n,r)}F \right) - \left( \Q_rF, \Q_rF \right)\right]<\eps_r,$$
где
$$ \Q_r = \frac 1 r\sum_{i=0}^{r-1} \tt^i.$$
Обозначая 
 $$P_r=\frac 1 {M(r)}\sum_{m=1}^{M(r)} \tt^{p(m,r)},$$
получим
$$\left\| \P_rF -\Q_rF \right\|^2  \ \leq \ \frac {\|F\|^2} {M(r)} + \eps_r.$$
Пусть   $M(r)$ и $p(m,r)$ выбраны  таким образом, что правая часть последнего
неравенства не превосходит величины, сравнимой с $2^{-r}$ (как выбрать, скажем позже).  
В этом случае  говорим, что векторы
$\P_{M(r)}F$  и  $\Q_rF $ очень близки, и пишем $\P_{M(r)}F\approx_r\Q_rF $.

Из сказанного выше и равенства
$$ I\otimes \T^r  +  \T^r\otimes I= (r+1)\Q_{r+1} -  r\Q_r  -  r \tt \Q_r  + (r-1)\tt \Q_{r-1}$$
получим
$$ (I\otimes \T^r  +  \T^r\otimes I)F\approx_r ((r+1)\P_{M(r+1)} -  r\P_{M(r)}  - $$
$$ - r \tt \P_{M(r)} \Q_r  + (r-1)\tt \P_{r-1})F:=\V_r F.\eqno (10)$$

Взяв в качестве $F$ вектор $f\otimes f$, получим
$$\T^{r} f\otimes f + f\otimes \T^{r}f\approx_r \V_r (f\otimes f).$$

Выше мы определили оператор $\V_r$, заметим, что его норма равна $4r$, поэтому очень близкие векторы 
он переводит в очень близкие. 
Подставляя в (10) вектор  $\T^{r} f\otimes f + f\otimes \T^{r}f$ вместо  $F$ и   $r+s$  вместо $r$,
получим
$$(\T^{r+s} \otimes I + I\otimes \T^{r+s})(\T^{r} f\otimes f + f\otimes \T^{r}f)\approx_{r+s}$$
 $$\V_{r+s}(\T^{r} f\otimes f + f\otimes \T^{r}f)\approx_{r}\V_{r+s}\V_r (f\otimes f).$$ 
Имеем
$$\tt^{r}(\T^{s} f\otimes f + f\otimes \T^{s}f)\approx_{r}\V_{r+s}\V_r (f\otimes f) - 
(\T^{2r+s} f\otimes f + f\otimes \T^{2r+s}f).$$

Справа стоит разность двух векторов, первый из которых принадлежит циклическому пространству $C_{f\otimes f}$,
а второй удален от  $C_{f\otimes f}$ на  малое расстояние, сравнимое с    $2^{-2r-s}$.
Получается, что для любого фиксированного $s>0$  вектор 
$\T^{s} f\otimes f  +  f \otimes  \T^{s} f$ отстоит от  $C_{f\otimes f}$
на расстоянии  не больше величины порядка $2^{-r }$, следовательно, расстояние равно нулю.

Осталось сказать, как выбраются  $M(r)$ и $p(m,r)$, $1\leq m\leq M(r).$

    Пусть монотонная последовательность $r_j$    принимает $2^{3r}$ раз подряд  значение $r$,
начиная с номера $j_r$. Положим $M(r)=2^{2r}$,   $p(m,r)= h_{j_r +2^{2r}m}$, где     $m=1,2, \dots, 2^r$.
Непосредственно проверяется, что  $\eps_r$ стремятся к 0 значительно быстрее, чем $2^{-r }$.
Таким образом, для данной конструкции имеем
$$\left\| \P_rF -\Q_rF \right\|  \ \leq 2^{-r}.$$

Для некоторых лестничных конструкций  верно, что их симметрические тензорные степени имеют простой спектр \cite{07}.
Гипотеза о том, что тем же свойством обладают лестничные конструции с логарифмическим ростом последовательности $r_j$ 
представляется весьма правдоподобной.

Читатель  может заметить, что логарифмический рост мы выбрали, чтобы избавиться от вычислений. На самом деле его можно
заменить на степенной.   Например, при $r_j^4\sim j$ аналогичные рассуждения также 
приводят к тому, что  векторы $\T^sf\otimes f + f\otimes \T^sf$
принадлежат циклическому пространству $C_{f\otimes f}$.  

Интересно узнать, какова кратность спектра  $\T\otimes\T$ (и других тензорных степеней)
 для лестничной конструкции $T$
в случаях  $r_j=[\sqrt j]$, $r_j=j$ и $r_j=j^2$?

Замечание. Известны преобразования  $R $ ранга один,  у которых $R\times R$ имеет  максимальную  
спектральную кратность, равную $2^n$. Они 
  строятся  следующим образом.  Находят $n$  преобразований, у которых кратность
спектра декартового квадрата равна 2,   их декартово произведение $R$ имеет ранг один, а   
 максимальная кратность спектра $R\times R$  равна $2^n$. 

\section{Бесконечные  преобразования, самоподобные конструкции}

Следуя \cite{2013c}, рассмотрим простейшую самоподобную конструкцию $T$, 
заданную параметрами $\bar s_j=h_j(0,1).$
Так как $$\T^{h_j}\to_w \frac 1 2 I,$$
 спектр  этой конструкции сингулярен, причем сверточные степени спектральной меры взвимно сингулярны.
Заметим, что 
$$\T^3\cong \T\oplus \T\oplus \T.$$
Степень $T^3$ не эргодична, пространство $X$ распадается на три инвариантных множества, на каждом из которых $T^3$
подобно $T$. Именно в этом смысле мы называем преобразование самоподобным.
По этой причине спектр $\T\otimes \T^3$  является сингулярным.
Будет ли спектр произведения $\T\otimes \T^2$ сингулярным? Этот  вопрос  оставим для изучения.

Рассмотрим другую конструкцию с параметрами  $\bar s_j=h_j(0,3).$
Произведение 
$T\otimes T^2$ диссипативно и по этой причине  имеет счетнократный лебеговский спектр. Произведение 
$T\otimes T^5$ консервативно и имеет сингулярный спектр.  
Похожие, но не самоподобные примеры приведены  в \cite{2019}.  Простейшие среди них имеют надстройки вида 
 $\bar s_j=(0, s_j),  \ s_j>>h_j.$

В работе  \cite{KuR} рассмотрена конструкция с  параметрами   
$$\bar s_j=h_j(1,5),  \ h_j=8^j.$$
Для нее была описана полугруппа слабых пределов степеней.  Она состоит из нулевого оператора  
и операторов вида
$2^{-m}\T^s$, $m=0,1,2,\dots$. В силу результатов \cite{2019}  эта конструкция обладает 
тривиальным централизатором.

\bf Асимметрия прошлого и будущего. \rm  У бесконечных самоподобных конструкций легко обнаружить свойство асимметрии.
В случае конечной меры  это требовало значительно  больше  изобретательности, см., например,  \cite{Os71}.

\bf Теорема 11.1.  \it Преобразование   с   параметрами   
$\bar s_j=h_j(0,1,2)$ неизоморфно своему обратному.\rm

Доказательство. Легко проверяется, что для   $\xi_i$-измеримого множества  $A$ конечной меры 
 для всех   $j>i$  выполнено:
$$\mu(A\cap T^{h_j}A\cap T^{3h_j}A)= \ \frac {\mu(A)} 3 ,$$
$$\mu(A\cap T^{2h_j}A\cap T^{3h_j}A)=0,$$
откуда вытекает
$$\mu(A\cap T^{-h_j}A\cap T^{-3h_j}A)=0.$$

Если $T$   и   $S$ сопряжены, то 
$$\mu(A\cap S^{h_j}A\cap S^{3h_j}A)\ \to \ \frac {\mu(A)} 3,$$
А если 
$T^{-1}$  и   $S$ сопряжены, то 
$$\mu(A\cap S^{h_j}A\cap S^{3h_j}A)\ \to 0,$$
поэтому $T$   и  $T^{-1}$ не могут иметь общее сопряжение $S$,
значит они не изоморфны.
 
\bf Самоподобные потоки.  \rm   Поток ранга один определяется похожим образом.  Башня на этапе $j$ 
отождествляется с прямоугольником высоты $h_j$, который разрезается по вертикали на $r_j>1$ колонн одинаковой ширины, 
над ними надстраиваются прямоугольники с высотами $s_j(i)\in \R$. Точка движется вертикально с постоянной скоростью,
достигнув верха надстроенной колонны с номером $i$, она  оказывается в основании колонны с номером $i+1$.
Мы мысленно складываем надстроенные колонны в одну, которая теперь есть прямоугольник этапа $j+1$ и т.д.(более подробное
описание см. в \cite{LR}).  Поток $T_t$  с надстройками 
$\bar s_j=h_j(0,\alpha)$  подобен потоку $T_{(2+\alpha) t}$  в том смысле, что 
они сопряжены в группе несингулярных преобразований. 
Гауссовские надстройки над таким потоком  пополняют коллекцию самоподобных потоков. 

Отметим, что  бесконечные преобразования ранга один чаконовского типа нашли интересное приложение 
к  проблеме изоморфизма пуассоновских надстроек  с одинаковым спектром \cite{RR}. При изменении  плотности меры,
в этом случае  пуассоновская надстройка меняется на неизоморфную, но с тем же  спектром. 
 В этом ключе было бы интересно изучить
надстройки над самоподобными потоками, предложенными выше.   Цель -- получить пуассоновский поток $P_t$ такой, что 
он не изоморфен $P_{\alpha t}$, хотя  у них одинаковый сингулярный спектр.  

Используя модификации конструкций из  работы \cite{2019}, можно построить  пуассоновский поток, 
удовлетворяющий условиям следующей теоремы.

\bf Теорема 11.2.  \it Для заданных непересекающихся счетных множеств $C,$ $C'$ найдется слабо перемешивающий поток 
$T_t$ такой, что спектр произведений $\T_1\otimes \T_c$   сингулярный при $c\in C$, $c>1$,
но имеет лебеговскую компоненту  при  $c\in C'$, $ c>1$.\rm   

Спектральная мера $\sigma$ такого потока обладает свойством:
произведение  $\sigma\times\sigma$ имеет как сингулярные так и несингулярные  проекции на диагональ в $\R^2$ в зависимости
от угла проекции.  Напрашивается такая терминология:  $\sigma\times\sigma$  прозрачна вдоль одних направлений и 
дает полную тень вдоль других.

\bf Перемешивающие конструкции. \rm При выборе параметров 
$$h_j<<s_j(1)<<s_j(2)<<\dots << s_j(r_j-1)<< s_j(r_j)$$
соответствующая конструкция обладает свойством перемешивания. 
 Действительно, пусть $f$ -- индикатор, например,    отрезка $E_1$. Он является циклическим вектором,
  а корреляции удовлетворяют неравенству:
$$|(\T^n f, f)| \leq \frac 1{r_j}, \ \   h_j<n\leq h_{j+1}.$$
Получаем  свойство перемешивания:    $\T^n\to 0, \ n\to\infty.$
Таким образом,  см. \cite{mmo14}, можно добится быстрых убываний  корреляций:
$$|(\T^n f, f)|< C n^{-\frac 1 2 +\eps}.$$
Свертка $\sigma\ast \sigma$ спектральной меры  $\sigma$ такого преобразования является лебеговской. Это
неудивительно, так как $T\times T$ -- диссипативное преобразование.  Действительно,  для сидоновких конструкций из условия 
$\sum_{j} \frac 1 {r_j}<\infty$  вытекает, что
точки из $E_1\times E_1$ со временем перестают возвращаться в $E_1\times E_1$.
И мы получаем, что $X\times X=\uu_{i\in\Z} (T\times T)^i Y$ для измеримого множества  $Y$ бесконечной меры.

Интересно узнать, какой спектральный тип меры $\sigma$ может быть у таких преобразований $T$. 
Если найти бесконечное преобразование ранга один с лебеговским спектром, то это решает известную задачу Банаха
о преобразовании с простым лебеговским спектром, упомянутую Уламом в  \cite{U}.

\bf Задача Тувено. \rm    Мы видим, что для сидоновских конструкций их декартовы квадраты -- диссипативные 
преобразования с бесконечным блуждающим множеством, поэтому они изоморфны. В связи с задачей  Ж.-П. Тувено
(см. \cite{2018}) возникает вопрос: какие инварианты могут различать сидоновские конструкции? Приведем пример.
 
Пусть $F_j\subset \N$ -- последовательность конечных множеств. Скажем, что преобразование $T$ бесконечного пространства 
с мерой принадлежит классу $\alpha\in [0,1]$  если для любого множества $A$ конечной меры
$$\limsup_j\mu\left(\cup_{n\in F_j}T^nA\ | \ A\right)=\alpha.$$
Очевидно, что $\alpha$ является инвариантом. Можно найти последовательность таких целочисленных интервалов 
 $F_j$, что  для любого 
$\alpha\in [0,1]$ найдется конструкция класса $\alpha$ с диссипативным декартовым квадратом.

В случае преобразований ранга один на вероятностном пространстве до сих пор не известно, 
влечет ли изоморфизм декартовх степеней
преобразований изоморфизм самих преобразований.

\section{ Заключительные замечания  и вопросы}

С преобразованиями ранга один связано большое количество открытых вопросов,
часть из них сформулирована ниже.

\bf 1. Общие вопросы о преобразованиях ранга один. \rm 

1.1. Дж. Кинг: верно ли, что любое эргодическое самоприсоединение $\nu$ преобразования $T$  ранга один  
есть предел внедиагональных мер: $\Delta^{k(j)} \to \nu$ ? Иначе говоря, верно ли, что  неразложимый марковский оператор 
$P$, коммутирующий с $\T$, лежит в слабом замыкании степеней $\T$?

 Мы знаем, что  для эргодического самоприсоединения $\nu$  найдется последовательность ${k(j)}$ такая, что
$\Delta^{k(j)} \to \eta\geq \frac 1 2 \nu.$ 
Гипотеза:  коэффициент $\frac 1 2 $ неулучшаем. 

1.2. Ж.-П. Тувено:  обладает ли  свойством  MSJ  нежесткое преобразование ранга один, у которого нет факторов?

1.3. Э. Руа:  верно ли, что ли   для бесконечных преобразований  ранга один, что их 
централизатор лежит в слабом замыкании степеней преобразования?

\bf 2.   Слабое кратное перемешивание WM(n). \rm  В теории  присоединений  давно открыт вопрос, 
связанный с проблемой о  кратном перемешивании: 
 \it 
верно ли, что слабо перемешивающее преобразование с нулевой энтропией 
не допускает нетривильные самоприсоединения с попарной независимостью? \rm

Хотя  имеется значительный  прогресс  (см. \cite{Host}), 
 для неперемешивающих действий ранга один вопрос  остается открытым.  
Ответ на  него  положителен для преобразования ранга один  со свойством  WMix(3) 
  слабого кратного  перемешивания   порядка 3. 
Напомним, что оно означает следующее: если  для  последовательностей 
$k_i, m_i, n_i$    имеет место сходимость
$$\mu(A\cap T^{k_i} B \cap T^{m_i} C\cap T^{n_i} D)\to \mu(A)\mu(B)\mu(C)\mu(D),$$
для всех наборов  множеств $A,B,C,D\in\B$, среди которых хотя бы одно совпадает с $X$,
то сходимость выполняется   для произвольных наборов $A,B,C,D\in\B$.
Свойство WM(n) определяется аналогично.

\bf Теорема 12. \it Для преобразований ранга один для любого фиксированного $k>0$ 
свойство WMix(2k+1) влечет за собой 
слабое кратное перемешивание WMix(m) для всех $m>1$ и отсутствие нетривиальных
самоприсоединений с попарной независимостью. \rm 

Этот факт  вытекает из результатов  \cite{R92},  \cite{97}. 

{\bf 3.   Гомоклинические группы преобразований ранга один. } 
   Гомоклиническая группа $ H (T)$  автоморфизма  $Т$  введена  М.И. Гординым:
$$ H (T) = \{S\in Aut(\mu): T^{- n} ST^n \to I, \ n \to \infty \}, $$
здесь подразумевается  сильная операторная сходимость к $I$.
Оказались содержательными некоторые  обобщения этого понятия. В частности, они представляют интерес в связи 
с общим вопросом:
насколько действия ранга один отличаются от гауссовских и пуассоновских действий (см. \cite{N}). Дело в том, что
последние обладают обширными гомоклиническими группами \cite{MMO19}, а действия ранга один в частных случаях
имеют  тривиальную гомоклиническую группы.  Дадим определения.

Слабо гомоклиническая  группа: 
$$ WH (T) = \{ S\in Aut(\mu): \frac {1} {N} \sum_{i = 1}^{N}T^{- i} ST^i \to I, 
\ N \to \infty \}. $$

Группа $wH(\{T_g\})$   действия $\{T_g\}$ (ее тоже можно назвать слабо гомоклинической) состоит из 
таких преобразований $S$, что  для всех последовательностей    
$\T_{g_i} \to_w \Theta$ выполняется   $\T^{- 1}_{g_i} \hat S\T_{g_i} \to I.$
 
Бесконечному множеству  $P\subseteq \Z$ и преобразованию  $T$ сопоставляется группа
$$ H_P (T) = \{S\in Aut(\mu) \ : \ \T^{- n} \hat S\T^n \to I, \ n \in P, \ n \to \infty \}. $$

Заметим, что имеют место включения: 

$H(T)\subseteq wH(\{T^n\}) \subseteq WH (T),$ \  $H(T)\subseteq  H_P (T).$

Для перемешивающих преобразований $T$ ранга один группа $WH (T)$ тривиальна \cite{MMO19}, но ничего не известно
про группы $H_P(T)$ для редких  множеств $P$.  Если группа $WH$  преобразования ранга один  эргодична,
то такое преобразование жесткое \cite{MMO19}. 

\it Какими бывают группы   $WH$ и $H_P$ для жестких слабо 
перемешивающих преобразований ранга один, для типичных преобразований?

\it Может ли  преобразование ранга один обладать  
эргодической слабо гомоклинической группой  $WH$?\rm

Отрицательный ответ подтверждает гипотезу о том, что пуассоновские надстройки не обладают  рангом один.
\\
\\
\bf 4. Конструкции в большим разбросом параметров. \rm 
Методы, обсуждавшиеся в этой статье,   не применимы к обширному, можно сказать, подавляющему множеству  конструкций 
ранга один. По ассоциации  с преобразованием Паскаля, введенным А.М.Вершиком в
\cite{V}, можно предложить  для изучения  спектральных и метрических свойств биномиальное преобразование  ранга один,
заданное параметрами 
 $r_j=j+1,  \ \ \ s_j(i)=C^i_j. $
   Обдумывая этот пример,  можно обнаружить изобилие конструкций, для которых   вопрос о непрерывности 
спектра, столь легкий для ограниченных конструкций,  является серьезной проблемой.  Вероятно, 
полезно расcмотреть   произведения вида
$$Q_J=\prod_{j\in J} P_j,\ \ P_j=\frac 1 {j+1} \sum_{i=0}^{j} T^{s_j(i)}$$
для больших конечных множеств $J$ (тут возникает ассоциация с произведениями Рисса).  Если большие  суммы 
из случайных биномиальных кэффициентов хорошо распределяются,
 возникает ситуация  $$  T^{-\sum_{j\in J}h_j}\approx_w Q_J\approx_w\Theta,$$
 что даст желаемую непрерывность спектра. Сложность представляет  последнее приближение. Первое 
при разумных ограничениях на $J$ является очевидным.

\bf 5. Несколько  задач о конструкциях ранга один. \rm

 5.1.  Для бесконечной лестничной конструкции  при выполнении условия  $\frac {r_j}{h_j}\to 0$
и ${r_j}\to \infty$  имеет место свойство перемешивания  \cite{stair}.
Обладает ли перемешиванием лестничная конструкция  для  $r_j$  при   $r_j=h_j$?  
Здесь возникают нетривиальные вопросы о том,  как квадратичные  подмножества $\N$  пересекаются со своими 
сдвигами. 

5.2.  Обладают ли простым спектром симметрические степени   ограниченных слабо перемешивающих конструкций? 
 
5.3. Верно ли, что 
что при $r_j\sim j^\alpha$, $\alpha > 0$, конструкции с надстройками $(0,s_j,0)$ являются простыми,
т.е. эргодические самоприсоединения порядка 2 суть $\mu\times\mu$ и меры вида $\Delta_S$?

  При $r_j=j$ слабое замыкание содержит множество  полиномиальных пределов, что
 говорит в пользу  минимальности спектра.

 При $r_j=j^n, n> 1,$  очевидны только пределы виды $a I + (1-a)\Theta$,  вопрос о существовании  пределов 
вида $a I + b\T^k+ \dots$ может оказаться трудной задачей.

5.4. Пусть  $p(j)$ обозначает  
 $j$-тое простое число. Глубокие факты  теории чисел приводят к тому, что 
конструкция $T$ с надстройками $(0,p(j),0)$
обладает  свойством дизъюнктности сверточных степеней  спектральной меры. 
  В  \cite{2013c} указаны нетривиальные слабые пределы как следствия новых  фактов о простых числах
\cite{Zh}.  Они обеспечивают свойство слабого перемешивания.  Из классических результатов
И.М. Виноградова следует, что для слабо перемешивающего преобразования $T$ имеет место сильная сходимость
$$\frac 1 {N} \sum_{j=1}^{N} T^{p(j)}\ \to \Theta.$$
Из неразложимости $\Theta$ в марковском централизаторе оператора $\T$, что эквивалентно слабому перемешиванию,
  для большинства простых 
чисел ${p}$ выполняется $T^{p}\approx_w\Theta$,  откуда для большинства $j$ получим
 $$T^{-h_j}\approx_w  \frac 1 2 I +\frac 1 3 T^{p(j)} +\frac 1 9 T^{p(j+1)}+\dots 
\approx_w\frac 1 2 I + \frac 1 2 \Theta.$$
А это влечет за собой  дизъюнктность сверток.
 Достаточен ли запас  слабых пределов  теоретико-числового происхождения  для того,
чтобы установить простоту спектра симметрических тензорных степеней оператора $\T$?

\bf 6. Полугруппы слабых пределов жестких преобразований. \rm 
Пусть $T$  -- преобразование,  $WLim(T)$ -- слабое замыкание действия  
$\{\T^n : n\in \Z\}$,  обозначим через $ULim(T)$ максимальную унитарную группу  в $WLim(T)$. 
Мы знаем примеры  неперемешивающих систем, для которых   $WLim(T)/ULim(T)$ имеет явное описание,
например, состоит из классов  $P^n$,  $n>0$,  где $P=\frac 1 2 (I+\T)$ или $P=\frac 1 2 (I+\Theta)$.
Однако, для жестких  слабо перемешивающих преобразований ранга один подобных описаний нет.
У конструкции дель Джунко-Рудольфа полугруппа  $WLim(T)$ содержит операторы $aI+(1-a)\T$  и $\Theta$. 
 Дают ли всевозможные их произведения  полное описание полугруппы $WLim(T)/ULim(T)$?

Хороший кандидат для изучения  $WLim(T)/ULim(T)$ -- жесткое преобразование ранга один, задаваемое
параметрами 
$$\bar s_j=(1,0,1,0, \dots,1,0,1,0), \ \ r_j\to \infty.$$

Другим кандидатом для исследования может служить \it $n!$-конструкция\rm, заданная параметрами
 $h_j=j!$, $r_j=j$, $s_j(i)=(j-1)!$ при $j>1$.
  Слабое замыкание такого бесконечного 
действия $T$ содержит операторы вида $a\hat S$, где $a\in [0,1]$ и   $\hat S\in ULim(T)$. Есть ли другие?

Гипотетическому ортогональному оператору  $V$, у которого $V^{n_i}\to I$ 
для некоторой последовательности ${n_i}\to\infty$, а 
слабое замыкание состоит из некоторой группы ортогональных операторов и
нулевого оператора, будет отвечать гауссовский  автоморфизм $G$ со слабым замыканием $WLim(G)=ULim(G)\cup \{\Theta\}$,
причем $ULim(G)$ -- континуальная группа.
Может ли таким же свойством обладать жесткое преобразование ранга один?
  

\end{document}